\documentclass{article}

\usepackage[top=30truemm,bottom=30truemm,left=25truemm,right=25truemm]{geometry}

\usepackage[utf8]{inputenc} % allow utf-8 input
\usepackage[T1]{fontenc}    % use 8-bit T1 fonts
\usepackage{hyperref}       % hyperlinks
\usepackage{url}            % simple URL typesetting
\usepackage{booktabs}       % professional-quality tables
\usepackage{amsfonts}       % blackboard math symbols
\usepackage{nicefrac}       % compact symbols for 1/2, etc.
\usepackage{fancyhdr}       % header
\usepackage{graphicx}       % graphics
\usepackage{dcolumn}% Align table columns on decimal point
\usepackage{bm}% bold math
\usepackage{comment,here,caption,subcaption,multirow}
\usepackage{amsmath,amsthm,amssymb}
\usepackage[all]{xy}

\usepackage{color}
\newcommand{\eref}[1]{(\ref{#1})}
\newcommand{\fref}[1]{Fig.~\ref{#1}}
\newcommand{\tref}[1]{Tab.~\ref{#1}}

\newtheorem{lemma}{Lemma}
\newtheorem{theorem}{Theorem}

\newtheorem*{remark*}{Remark}

\title{Koopman Analysis of the Singularly-Perturbed van der Pol Oscillator
\thanks{The work was partially supported by JSPS KAKENHI Grant Number 23H01434 and JST Moonshot R\&D Grant Number JP-MJMS2284.}
}

\author{
  Natsuki Katayama\thanks{Electronic mail : \texttt{n-katayama@dove.kuee.kyoto-u.ac.jp}}, Yoshihiko Susuki\thanks{Electronic mail : \texttt{susuki.yoshihiko.5c@kyoto-u.ac.jp}} \\
  Department of Electrical Engineering, Kyoto University 
}

\begin{document}
\maketitle

\begin{abstract}
The Koopman operator framework holds promise for spectral analysis of nonlinear dynamical systems based on linear operators.
Eigenvalues and eigenfunctions of the Koopman operator, so-called Koopman eigenvalues and Koopman eigenfunctions, respectively, mirror global properties of the system's flow.
In this paper we perform the Koopman analysis of the singularly-perturbed van der Pol system.
First, we show the spectral signature depending on singular perturbation: how two Koopman {principal} eigenvalues are ordered and what distinct shapes emerge in their associated Koopman eigenfunctions. 
Second, we discuss the singular limit of the Koopman operator, which is derived through the concatenation of Koopman operators for the fast and slow subsystems. 
From the spectral properties of the Koopman operator for the {singularly}-perturbed system and the singular limit,
we suggest that the Koopman eigenfunctions inherit geometric properties of the singularly-perturbed system. 
These results are applicable to general planar singularly-perturbed systems with stable limit cycles.
\end{abstract}

\section{Introduction}
% Intro to KOT
The Koopman operator framework has gained successful results as a novel way of analyzing nonlinear dynamics, see, e.g., Refs\cite{analysis_of_fluid_flowd_KO,brunton2021modern,otto2021koopman}. 
The Koopman operator is a linear composition operator defined for a large class of nonlinear dynamical systems (flows in continuous time and maps in discrete time) \cite{lasota1998chaos}. 
Even if a target system is finite-dimensional and nonlinear, the Koopman operator is infinite-dimensional but linear that completely captures the information on the target nonlinear system.  
Specifically, spectral properties of the Koopman operator have been linked to dynamics emerging in the nonlinear system, see, e.g., Refs\cite{arnold1968ergodic,lasota1998chaos,mezic2005spectral,rowley2009spectral} and references in Ref\cite{koopman-springer}. 
The eigenvalues of the Koopman operator and associated eigenfunctions, referred to as Koopman eigenvalues and {eigenfunctions}, play an important role in analyzing global flows described by systems with stable limit cycles  \cite{use_of_Fourier_averages,isostable_reduction_of_periodic_orbits,phase-amplitude_reduction_nakao,global_computation_of_phase-amplitude_reduction}, which we call the Koopman analysis.

In this paper we perform the Koopman analysis of a singularly-perturbed (SP) system with a stable limit cycle.
{SP systems have been traditionally} studied for analysis of multiple time scale dynamics in the physical world, for example, relaxation oscillations \cite{mishchenko2013differential,izhikevich2000phase,desroches2012mixed,kuehn2015multiple}.
The most famous example of such SP systems for relaxation oscillations is the SP van der Pol equation (see, e.g., Refs\cite{mishchenko2013differential,izhikevich2000phase}), given as 
\begin{equation}
  \label{eq:vanderPol_fast}
  \left.
    \begin{alignedat}{4}
      x' & {}:={} &
     {{\rm d} x \over {\rm d} t}  &{}={}& 
     F(x,y) &{}:= x - \frac{x^3}{3} + y\\
     y' &{}:={} &
      {{\rm d} y \over {\rm d} t} &{}={}& \varepsilon G(x,y) &{}:= -\varepsilon x
    \end{alignedat}
  \right\},
\end{equation}
where $t \in\mathbb{R}$ is the time,  $x, y\in \mathbb{R}$ are the state variables, and $\varepsilon\,(\ll 1)$ is the small positive parameter. 
By changing the time scale $\tau := \varepsilon t$, the equivalent system is obtained as
\begin{equation}
  \label{eq:vanderPol_slow}
  \left.
    \begin{alignedat}{3}
    \varepsilon \dot{x} &{}:={}& \varepsilon {{\rm d} x \over {\rm d} \tau} &{}={} &F(x,y)\\
     \dot{y} &{}:={}&  {{\rm d} y \over {\rm d} \tau} &{}={} & G(x,y)
    \end{alignedat}
  \right\},
\end{equation}
where $t$ represents the {\em fast} time scale while $\tau$ represents the {\em slow} time scale.
{Eq. \eref{eq:vanderPol_fast}} discribes the fast dynamcis by the $x$ variable of the SP system in $t\ll 1/\varepsilon$, while {Eq. \eref{eq:vanderPol_slow}} does the slow dynamics mainly characterized by the $y$ variable in $\tau\gg \varepsilon$.
For general SP systems, including the forms of (\ref{eq:vanderPol_fast}) and (\ref{eq:vanderPol_slow}), their operator-theoretic analysis has been reported in the literature. 
For example, the authors of Ref\cite{isostable_isochron_and_koopmanSpectrum_fixedpoint} discussed a connection between Koopman eigenfunctions and geometric properties of {an} invariant manifold for a system with a stable equilibrium point. 
The authors of Ref\cite{eldering2018global} showed that the conventional local Fenichel theory \cite{fenichel1979geometric,singular_perturbation_geometric} can be globally extended by virtue of the Koopman operator theory in Ref\cite{linearization_in_the_large_of_nonlinear_systems_and_KoopmanOperatorSpectrum}. 
The authors of Refs\cite{osinga2010continuation,sherwood2010dissecting,mauroy2014global} studied the geometry of isochrons for slow-fast systems with limit cycles, which are mathematically defined as level sets of a particular Koopman eigenfunction. 
The above-existing research focuses on smooth dynamics described by the SP systems. 

\begin{figure}[h]
    \centering
    \begin{minipage}{1\textwidth}
        \centering
        \[
    \xymatrix@!C=120pt{
    {\left\{ \begin{alignedat}{3}\varepsilon {{\rm d} x \over {\rm d} \tau} &{}={}& \varepsilon \dot{x} &{}={} &F(x,y)\\{{\rm d} y \over {\rm d} \tau} &{}={}& \dot{y} &{}= & G(x,y)\end{alignedat} \right.}
    \ar[r]^*{t = \tau/\varepsilon} \ar[d]_*{\varepsilon = 0} & \ar[l] \ar[d]_*{\varepsilon = 0}  
    {\left\{ \begin{alignedat}{2}x' =  {{\rm d} x \over {\rm d} t} &{}={} &F(x,y)\\y' =  {{\rm d} y \over {\rm d} t} &{}= & \varepsilon G(x,y)\end{alignedat} \right. }
    \\ 
    {{\rm slow :}\left\{ \begin{aligned}0 &{}= F(x,y) \\ \dot{y} &{}= G(x,y)\end{aligned} \right.}
    & 
    {{\rm fast :}\left\{ \begin{aligned}    x' &{}= F(x,y) \\ y' &{}= 0\end{aligned} \right.}
   } 
        \]
        \caption{Relation among two singularly-perturbed systems and associated slow/fast subsystems.}
        \label{fig:ODE_diagram}
    \end{minipage}
\end{figure}

{The novelty of our Koopman analysis} is to characterize spectral signatures of singular perturbation in the van der Pol equation, including the singular limit \cite{kuehn2015multiple} as $\varepsilon\to 0$, where nonsmooth dynamics evolve over the whole state plane. 
At the singular limit, the asymptotic solutions of the SP system become non-smooth or discontinuous in terms of $\tau$.
As shown in \fref{fig:ODE_diagram}, the singular limit induces the decomposed two subsystems (\ref{eq:fast_subsystem}) and (\ref{eq:slow_subsystem}) from (\ref{eq:vanderPol_fast}) and (\ref{eq:vanderPol_slow}) as follows: 
\begin{equation}
  \label{eq:fast_subsystem}
  \left.
    \begin{aligned}
      x' & =F(x,y)\\
      y' & = 0
    \end{aligned}
  \right\},
\end{equation}
and
\begin{equation}
  \label{eq:slow_subsystem}
  \left.
    \begin{aligned}
      0 & =F(x,y)\\
      \dot{y} & = G(x,y)
    \end{aligned}
  \right\}.
\end{equation}
The {\em fast subsystem} \eqref{eq:fast_subsystem} is a form of Ordinary Differential Equation (ODE) and represents the fast evolution ${x}(t)$ on the boundary layer parameterized by a constant $y$,
while the {\em slow subsystem} \eqref{eq:slow_subsystem} is a form of semi-explicit Differential-Algebraic Equation (DAE) and represents the slow evolution $(\bar{x}(\tau),\bar{y}(\tau))$ on the critical manifold as
\begin{equation}
\label{eq:critical_manifold}
 W := \{(x,y)\in\mathbb{R}^2 ~|~F(x,y)=0\}. 
\end{equation}
Then, the discontinuous solutions at the singular limit can be decomposed into solutions described by the fast and slow subsystems. 
In this paper, we clarify the operator framework for the singular perturbation including the discontinuity.

The main contributions of this paper are two-fold. 
First, we show a spectral signature of the singular perturbation, that is, how spectra of the Koopman operators for the SP systems depend on the small parameter $\varepsilon$. 
Precisely, if the SP van der Pol system has a stable limit cycle with fundamental angular frequency $\omega^\varepsilon$ and negative Floquet exponent $\nu^\varepsilon$, then ${\rm i}\omega^\varepsilon$ and $\nu^\varepsilon$ coincide with the two Koopman {principal} eigenvalues ordered $\mathcal{O}(1),~\mathcal{O}(\varepsilon^{-1})$, respectively. 
The numerical results of the two Koopman {principal} eigenfunctions suggest that they exhibit distinct shapes on the state plane of the fast and slow variables: there exist scale differences depending on $\varepsilon$ between both of their derivatives with respect to (w.r.t.) $x$ (fast) and $y$ (slow).
Regarding the eigenfunction associated with the slow eigenvalue ${\rm i}\omega^\varepsilon$, 
the change in the $y$ direction significantly appears in a neighborhood of the critical manifold $W$, where the shape of the Koopman eigenfunction becomes steep or sharp depending on $\varepsilon$.

Second, we provide a formalism of Koopman operators for the SP van der Pol oscillator {systems} that include the case of the singular limit. 
For the fast subsystem (\ref{eq:fast_subsystem}) which is an ODE parameterized by $y$, the Koopman operator and its spectral properties are naturally introduced without any novelty. 
On the other hand, for the slow subsystem (\ref{eq:slow_subsystem}) which is a semi-explicit DAE, we show that the Koopman operator is well-defined on $W$, which is beyond the existing definition in Ref\cite{susuki2021koopman}.  
As its spectral property, it is also shown that one Koopman eigenvalue is regarded as a limit of ${\rm i}\omega^\varepsilon$ as $\varepsilon\to 0$. 
Then, by combining these Koopman operators, we construct the Koopman operator at the singular limit that is well-defined on the entire state plane. 
The formalism suggests that the spectral signature to the singular perturbation shown above is involved in the spectral property at the singular limit.
This suggests that the distinctive shapes of the Koopman eigenfunction capture the geometric properties of the SP flow.

The rest of the paper is organized as follows.
Section~\ref{sec:II} will provide a brief introduction to the SP van der Pol system and derive the flows for the slow and fast subsystems, and the singular limit. 
In Section III we will introduce the Koopman operators for the SP van der Pol system and discuss its spectral properties.
In Section IV we will introduce the Koopman operators for the slow and fast subsystems, and the singular limit, and we will derive their Koopman eigenfunctions. 
Conclusions are presented in Section~V with a brief summary and future directions.

%%%%%%%%%%%%%%%%%%%%%% Sec.II %%%%%%%%%%%%%%%%%%%%%%%%
\section{Geometry and flows of the Singularly-Perturbed van der Pol Oscillator}
\label{sec:II}
\begin{figure*}[t]
  \begin{tabular}{cc}
    \begin{minipage}[t]{0.5\hsize}
      \centering
      \includegraphics[keepaspectratio, scale=0.7]{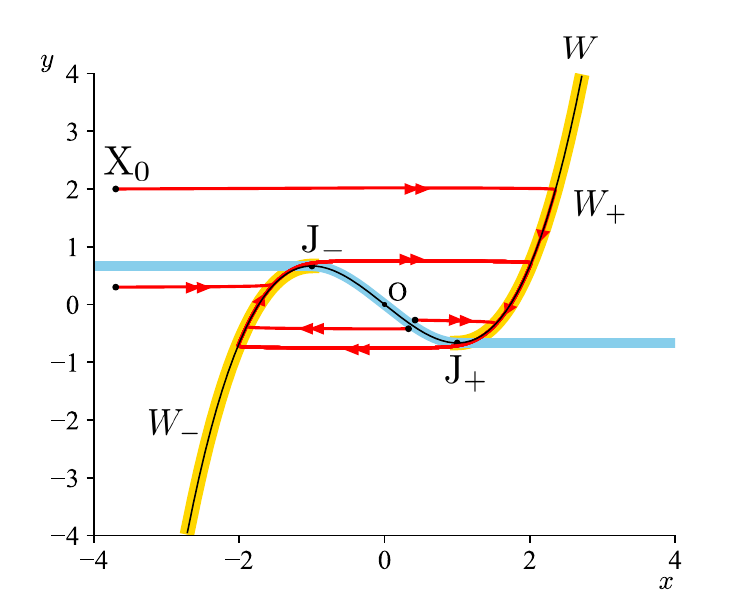}
      \caption{State plane of (\ref{eq:vanderPol_slow}) where $\varepsilon=0.01$. $W$ is represented by the black curve and $W_\mp$ are painted in yellow.}
      \label{fig:vander001}
    \end{minipage} &
    \begin{minipage}[t]{0.5\hsize}
      \centering
      \includegraphics[keepaspectratio, scale=0.7]{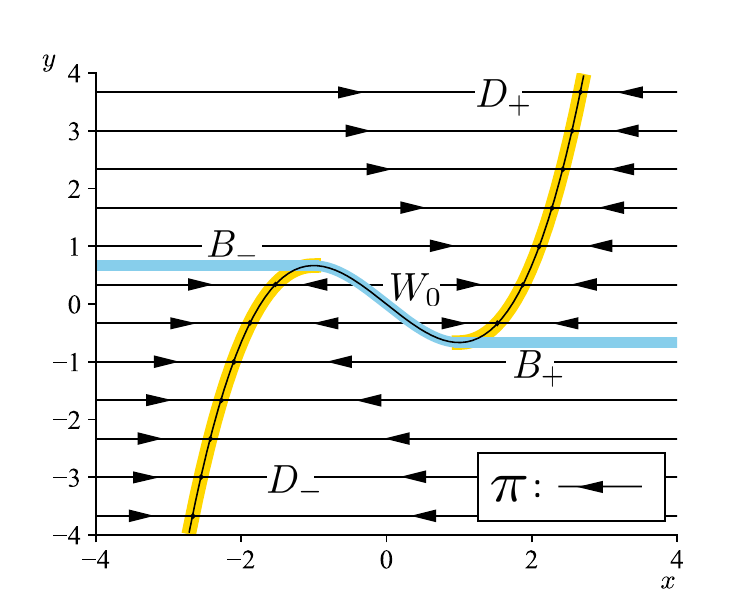}
      \caption{$B_-,~W_0,~B_+$ (forming the thick curve colored by skyblue), $D_-$ and $D_+$ separated by the curve $B_-\cup W_0 \cup B_+$, and projection $\pi$ (black lines with arrows).}
      \label{fig:projection}
    \end{minipage}
  \end{tabular}
\end{figure*}

This section briefly introduces the dynamical systems (\ref{eq:vanderPol_fast}) and \eref{eq:vanderPol_slow} of the SP van der Pol oscillator. 
We mainly analyze the $\tau$-scaled ODE (\ref{eq:vanderPol_slow}) in the rest of the paper. 

The state plane of the system (\ref{eq:vanderPol_slow}) is illustrated in \fref{fig:vander001}. 
For sufficiently small $\varepsilon$, the solution of (\ref{eq:vanderPol_slow}) from ${\rm X}_0$ moves fast in the $x$ direction, reaches in a neighborhood of the critical manifold $W$ as (\ref{eq:critical_manifold}),
and then moves slowly along $W$, that is, along the slow manifold \cite{kuehn2015multiple}.  
Finally, it converges to the stable limit cycle that represents the relaxation oscillation. 
To describe this from a geometrical viewpoint in \fref{fig:vander001}, 
let us introduce the three subsets of $W$ as 
$$
\left. \begin{aligned}
    W_- &:= W \cap \{ (x,y)\in \mathbb{R}^2 ~|~x<-1 \} \\
    W_0 &:= W \cap \{ (x,y)\in \mathbb{R}^2 ~|~-1\leq x\leq 1 \} \\
    W_+ &:= W \cap \{ (x,y)\in \mathbb{R}^2 ~|~x>1 \}
\end{aligned} \right\}.
$$
Also in \fref{fig:vander001}, the two points ${\rm J}_- = (- 1, + 2/3)$ and ${\rm J}_+ = (+ 1,- 2/3)$ are introduced and called the singular points (otherwise known as jump points \cite{kuehn2015multiple}). 
The solution described above is then geometrically interpreted as follows: the solution from ${\rm X}_0 =(-3.8,2)$ (as shown in \fref{fig:vander001})
  {
  (i) moves fast in the $x$ direction toward $W_+$ and enters the neighborhood of $W_+$;
  (ii) moves slow along $W_+$ toward ${\rm J}_+$ and enters the neighborhood of ${\rm J}_+$; 
  (iii) leaves the neighborhood of ${\rm J}_+$, moves fast toward $W_-$, and enters the neighborhood of $W_-$; 
  (iv) moves slow along $W_-$ toward ${\rm J}_-$ and enters the neighborhood of ${\rm J}_-$;
  (v) leaves the neighborhood of ${\rm J}_-$, moves fast toward $W_+$, and enters the neighborhood of $W_+$;} and 
  (vi) repeats (ii)-(v).
Regarding (i), which the solution heads to $W_+$ or $W_-$, depends on the location of the initial state. 
The location is roughly determined with $W_0$ and the two half-lines $B_-$ and $B_+$, defined as
$$
\left. \begin{aligned}
B_- &= \left\{ \left(x,y \right)\in \mathbb{R}^2 ~\middle|~x < -1,~y=2/3 \right\} \\
B_+ &= \left\{ \left(x,y \right)\in \mathbb{R}^2 ~\middle|~x > 1,~y=-2/3 \right\}
\end{aligned} \right\}.
$$
The sets $W_-$, $W_+$ are depicted in \fref{fig:vander001}, while $W_0$ and $B_\mp$ are depicted in \fref{fig:projection}.
Below, we denote the semi-group of flows for the SP systems (\ref{eq:vanderPol_slow}) and (\ref{eq:vanderPol_fast}) as $\{ \boldsymbol{S}_{{\rm slow},\tau}^\varepsilon :\mathbb{R}^2\backslash \{\boldsymbol{0} \}\to \mathbb{R}^2\backslash \{\boldsymbol{0} \} ,~\tau \geq 0\}$ and $\{ \boldsymbol{S}_{{\rm fast}, t}^\varepsilon:\mathbb{R}^2\backslash \{\boldsymbol{0} \}\to \mathbb{R}^2\backslash \{\boldsymbol{0} \},~t \geq 0\}$, satisfying 
\begin{equation}
\boldsymbol{S}_{{\rm slow},\tau}^\varepsilon = \boldsymbol{S}_{{\rm fast},\tau/\varepsilon}^\varepsilon,  
\qquad
 \boldsymbol{S}_{{\rm slow},\varepsilon t}^\varepsilon = \boldsymbol{S}_{{\rm fast},t}^\varepsilon.
 \label{eqn:conversion_flow}
\end{equation}
This implies that the $\tau$-scaled flows $\boldsymbol{S}^\varepsilon_{{\rm slow},\tau}$ and $\boldsymbol{S}^\varepsilon_{{\rm fast},\tau/\varepsilon}$ generate common solutions in $\tau$ scale, namely, the slow time-scale.

Next, we introduce the flows of the fast subsystem (\ref{eq:fast_subsystem}) and the slow one (\ref{eq:slow_subsystem}). 
%%%%%% fast
The fast subsystem (\ref{eq:fast_subsystem}) is a 1-dimensional system with equilibrium points located at $W$, where $y$ is a constant, and equilibrium points on $W_\mp$ are attracting while equilibrium points on $W_0$ are repelling.
We denote the flow of the fast subsystem (\ref{eq:fast_subsystem}), named the {\em fast flow}, by $\boldsymbol{S}_{{\rm fast}, t}^0(\,\cdot\, ;y):\mathbb{R}\rightarrow \{(x,y)\in \mathbb{R}^2 ~|~ y = {\rm const.}\}~(t\geq 0)$, i.e., $\boldsymbol{S}_{{\rm fast}, t}^0 (x;y)$ is a solution to (\ref{eq:fast_subsystem}) with the initial state  $(x;y)$. 
Note that the fast flow parameterized as $t\geq 0$ is the map from the domain $\mathbb{R}$ to the range $\mathbb{R}^2$ including constant $y$. 
%%%%%% slow
The slow subsystem (\ref{eq:slow_subsystem}) represents the semi-explicit DAE on $W_\mp$. 
For introducing the state space of \eqref{eq:slow_subsystem}, we exclude $W_0$ by following Ref\cite{takens1976constrained}.
In the slow subsystem of \eref{eq:slow_subsystem}, any solutions can no longer be defined beyond some $\tau_{\rm s}$, which represents the time it takes to reach the singular points ${\rm J}_\mp$ from initial states.
In constrained differential system theory \cite{takens1976constrained,sastry1981jump}, when $\tau = \tau_{\rm s}$, the solutions jump discontinuously from the singular points to the points $(2,2/3)\in W_+$ and $(-2,-2/3)\in W_-$, called drop points.
By following this theory, the flow on $W_\mp := W_-\cup W_+$, named the {\em constrained flow}, is defined as $\boldsymbol{S}_{{\rm slow},\tau}^0 : W_\mp \to W_\mp~(\tau\geq 0)$.
Its precise definition is given in Appendix A.

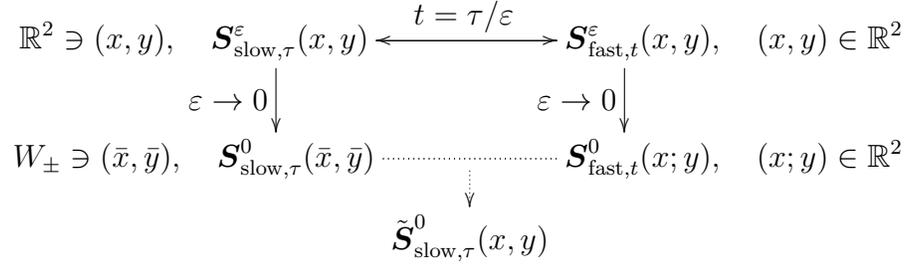
\begin{figure*}[tb]
    \centering
    \begin{minipage}{1\textwidth}
       {\large \[
    \xymatrix@C=0pt@R=9pt{
    {\mathbb{R}^2\ni (x,y),\quad \boldsymbol{S}_{{\rm slow},\tau}^\varepsilon (x,y)}
    \ar[rr]^*{t = \tau/\varepsilon} \ar@<6ex>[dd]_*{\varepsilon \to 0} & & \ar[ll] \ar@<-8ex>[dd]_*{\varepsilon \to 0}  %矢印の設定
    {\boldsymbol{S}_{{\rm fast},t}^\varepsilon (x,y), \quad (x,y) \in \mathbb{R}^2 }
    \\ & & \\ 
    {W_\pm \ni (\bar{x},\bar{y}),\quad \boldsymbol{S}_{{\rm slow},\tau}^0 (\bar{x},\bar{y})}
    \ar@{.}[rr] & \ar@{.>}[d]&\ar@{.}[ll] 
    {\boldsymbol{S}_{{\rm fast},t}^0 (x;y), \quad (x;y) \in \mathbb{R}^2}
    \\
    & {\tilde{\boldsymbol{S}}_{{\rm slow},\tau}^0 (x,y)} &
   }
        \]}
        \caption{{\color{black}Five flows derived from the fast and slow subsystems, and their limits which we analyze in this paper.}}
        \label{fig:flow_diagram}
    \end{minipage}
\end{figure*}
Finally, we define the flow of the singular limit. 
For this, several notations are introduced first. 
For the subsets of $W$, there exist graphs $\gamma_\mp,~\gamma_0$ such that $W_\mp,~W_0$ can be represented as 
$$
\left. \begin{aligned}
W_- &= \{ (x,y)\in\mathbb{R}^2~|~x=\gamma_-(y),~y<2/3 \} \\
W_0 &= \{ (x,y)\in\mathbb{R}^2~|~x=\gamma_0(y),~ {-2/3\leq y\leq 2/3} \} \\
W_+ &= \{ (x,y)\in\mathbb{R}^2~|~x=\gamma_+(y),~y>-2/3 \}
\end{aligned} \right\}.
$$
To characterize the asymptotics of the fast subsystem, we introduce the subdomains $D_\mp$ of the state plane  $\mathbb{R}^2$ as
$$
\left. \begin{aligned}
    D_- &= \{ (x,y)\in \mathbb{R}^2 ~|~ (x < \gamma_0 (y) \cap |y| \leq 2/3)\cup (y<-2/3) \} \\
    D_+ &= \{ (x,y)\in \mathbb{R}^2 ~|~ (x > \gamma_0 (y) \cap |y| \leq 2/3)\cup (y>2/3) \} \\
\end{aligned} \right\},
$$
which are separated by the curve $W_0$ and $B_\mp$ connected (see \fref{fig:projection}). 
If an initial state is in $D_+$ (or $D_-$), then the associated solution of the fast subsystem asymptotically converges to $W_+$ (or $W_-$) as $t\to \infty$. 
Here, as described in the introduction, the discontinuous solutions at the singular limit are decomposed into the fast and slow subsystems. 
For their decomposition in $\tau$ scale, the $\tau$-scaled flow $\boldsymbol{S}^\varepsilon_{{\rm fast}, \tau/\varepsilon}$ is regarded as a map as $\varepsilon\to 0$: precisely, it becomes the projection $\pi$ shown in \fref{fig:projection} from $\mathbb{R}^2\backslash W_0$ to $W_\mp$ satisfying 
\begin{equation}
\label{eq:projection}
\pi (x,y) = \left\{ \begin{aligned}
    &(\gamma_+ (y), y) ,\quad (x,y)\in D_+,\\
    &(\gamma_- (y), y) ,\quad (x,y)\in D_-.
\end{aligned}\right. 
\end{equation}
Thus, as the combination of $\pi$ and the constrained flow $\boldsymbol{S}^0_{{\rm slow},\tau}$, the flow at the singular limit is defined as 
$\tilde{\boldsymbol{S}}_{{\rm slow}, \tau}^0:\mathbb{R}^2\backslash W_0 \rightarrow \mathbb{R}^2\backslash W_0$, satisfying 
\begin{equation}
    \label{eq:singularflow}
    \tilde{\boldsymbol{S}}_{{\rm slow}, \tau}^0 (x,y) := \left\{ \begin{alignedat}{2} & (x,y), & \tau = 0,\\
        & \boldsymbol{S}_{{\rm slow},\tau}^0 \circ \pi (x,y), \quad & \tau > 0.
    \end{alignedat}\right.
\end{equation}
The flow $\tilde{\boldsymbol{S}}_{{\rm slow},\tau}^0$ generates the discontinuous solutions in $\tau$ scale (slow time-scale). 
The relation among the five flows is summarized in \fref{fig:flow_diagram}.

\section{Koopman Operators of the Singularly-Perturbed Systems}
In this section, we introduce the Koopman operators for the SP van der Pol system with $\varepsilon > 0$ and discuss their spectral properties. 

\subsection{Slow/fast Koopman operators}
First of all, we introduce the Koopman operator for the flow $\boldsymbol{S}^\varepsilon_{{\rm slow}, \tau}$ of the SP system \eqref{eq:vanderPol_slow} and its eigenfunction. 
The Koopman operator is a composition operator of the %system's 
flow and observables that are functions from a dense set 
$\mathbb{R}^2\backslash\{\boldsymbol{0}\}$ of the state plane to a scalar. 
Following Ref\cite{spectrum-of-the-koopman-operator} to derive a nice spectral property as explained below, we take the space $\mathcal{F}$ of observables as consisting of all analytic functions on $\mathbb{R}^2\backslash\{\boldsymbol{0}\}$.  
The semi-group of Koopman operators, $\{U_{{\rm slow},\tau}^\varepsilon :\mathcal{F} \rightarrow \mathcal{F}, \tau\geq 0\}$, 
is then defined as follows: for $f\in\mathcal{F}$,
\begin{align}
  \label{eq:koopmndefine}
  (U_{{\rm slow},\tau}^\varepsilon f) (x,y) 
  &:= (f\circ \boldsymbol{S}_{{\rm slow}, \tau}^\varepsilon ) (x,y) \\
  &= f(\boldsymbol{S}_{{\rm slow},\tau}^\varepsilon (x,y)). \nonumber
\end{align}
Note that at $\tau =0$, $U_{{\rm slow},\tau}^\varepsilon$ becomes an identity operator: $U_{{\rm slow},0}^\varepsilon f = f$ for all $f$.
An eigenvalue $\lambda\in\mathbb{C}$ of the operators' group and associated eigenfunction  
$\phi_\lambda^\varepsilon \in \mathcal{F}\backslash\{ 0 \}$ can be defined as 
\begin{equation}
  \label{eq:KEFdef}
  (U^\varepsilon_{{\rm slow},\tau} \phi_\lambda^\varepsilon)(x,y) = {\rm e}^{\lambda \tau} \phi_\lambda^\varepsilon (x,y). 
\end{equation}
In literature (see, e.g., Refs\cite{koopman-springer,spectrum-of-the-koopman-operator}), $\lambda$ is called the Koopman eigenvalue, and  $\phi_\lambda^\varepsilon$ the associated Koopman eigenfunction. 
{For an analytic system with a stable limit cycle such as the SP van der Pol system, it is stated from Ref\cite{spectrum-of-the-koopman-operator} that the operator's group has a countably infinite number of pairs of Koopman eigenvalues and eigenfunctions. }
In what follows, 
we focus on the so-called {principal} pairs that are utilized to
characterize the system's properties \cite{spectrum-of-the-koopman-operator}. 
The Koopman principal eigenvalues are ${\rm i}\omega^\varepsilon$ and $\nu^\varepsilon$ (i: imaginary unit),  
which are associated with the fundamental angular frequency $\omega^\varepsilon$ and the Floquet exponent $\nu^\varepsilon$ of the limit cycle, respectively. 
The associated {principal} eigenfunctions $\phi_{{\rm i}\omega}^\varepsilon$ and $\phi_\nu^\varepsilon$ will be visualized in the next subsection. 
The Koopman {principal} eigenfunctions are unique \cite{Existence_Kvalheim}, that is, analytic functions satisfying {\eqref{eq:KEFdef}} with $\lambda={\rm i}\omega^\varepsilon$ (or $\nu^\varepsilon$) are only $\phi_{{\rm i}\omega}^\varepsilon$ (or $\phi_\nu^\varepsilon$) and its constant multiples.
Such {principal} eigenfunctions inherit the geometric properties of the limit cycling system as follows:
For $\phi_{{\rm i}\omega}^\varepsilon$, one can represent $\phi_{{\rm i}\omega}^\varepsilon = {\rm e}^{{\rm i}\theta^\varepsilon}$ where $\theta^\varepsilon$ is the asymptotic phase of the limit cycle, which implies that $|\phi_{{\rm i}\omega}^\varepsilon|$ is constant and level sets of $\phi_{{\rm i}\omega}^\varepsilon$ correspond to the isochrons of the limit cycling system \cite{use_of_Fourier_averages}.
For $\phi_\nu^\varepsilon$, it is a real-valued function with zero value on the limit cycle, whose level sets correspond to the isostables  \cite{isostable_reduction_of_periodic_orbits,phase-amplitude_reduction_nakao}.

In the same manner as above, we can introduce the Koopman operator for the flow $\boldsymbol{S}^\varepsilon_{{\rm fast}, t}$ of the SP system {\eqref{eq:vanderPol_fast}}.
The semi-group of Koopman operators, $\{U_{{\rm fast},t}^\varepsilon :\mathcal{F} \rightarrow \mathcal{F}, t\geq0\}$, is defined as
\begin{align}
  \label{eq:koopmndefine_fast}
  (U_{{\rm fast},t}^\varepsilon f) (x,y) 
  &:= (f\circ \boldsymbol{S}_{{\rm fast},t}^\varepsilon ) (x,y), 
  \quad f\in\mathcal{F}.
\end{align}
The eigenvalue and eigenfunction of this $U_{{\rm fast},t}^\varepsilon$ are defined in the same manner as above and will appear in the next subsection. 
The Koopman analysis of the slow/fast subsystems \eqref{eq:slow_subsystem} and \eqref{eq:fast_subsystem}, and at the singular limit will be performed in Section~IV.

%%%%%
\subsection{Slow/fast Koopman eigenvalues}
Here and in the next subsection, we discuss the relationship between the spectrum of the Koopman operator and singular perturbation for the van der Pol system \eref{eq:vanderPol_slow}.
First, we clarify how the spectral relation of the two Koopman operators $U^\varepsilon_{{\rm slow},\tau}$ and $U^\varepsilon_{{\rm fast},t}$ for the slow and fast timescales are connected.
By using %$t=\tau/\varepsilon$ and 
\eqref{eqn:conversion_flow}, the two Koopman operators are connected in $t$ scale as: for $f\in\mathcal{F}$
\[
U_{{\rm fast},t}^\varepsilon f 
= f\circ\boldsymbol{S}^\varepsilon_{{\rm fast},t}
= f\circ\boldsymbol{S}_{{\rm slow}, \varepsilon t}^\varepsilon 
= U_{{\rm slow}, \varepsilon t}^\varepsilon f, %.
\]
and in $\tau$ scale as 
\[
U_{{\rm slow},\tau}^\varepsilon f
= U_{{\rm fast},\tau/\varepsilon}^\varepsilon f.
\]
Thus, \eqref{eq:KEFdef} is written in $t$ scale as
$$
\begin{aligned}
    {U}_{{\rm fast}, t}^\varepsilon \phi_\lambda^\varepsilon 
    = U_{{\rm slow}, \varepsilon t}^\varepsilon \phi_{\lambda}^\varepsilon 
    = {\rm e}^{\varepsilon \lambda t} \phi_{\lambda}^\varepsilon.
\end{aligned}
$$
This implies
if (\ref{eq:vanderPol_slow}) in $\tau$ scale has the Koopman eigenvalue $\lambda$, then (\ref{eq:vanderPol_fast}) in $t$ scale has $\varepsilon \lambda$. 
Their eigenfunctions are the same, that is, $\phi^\varepsilon_\lambda$ (see \tref{tab:timescale}).
\begin{table}[t]
  \caption{
  Time scale transformation of Koopman operators, Koopman eigenvalues, and Koopman eigenfunctions.
  }%
  \label{tab:timescale}
    \begin{center}
    \begin{tabular}{cccc}
    ~ & $\tau$-scale & $\longleftrightarrow$ & $t$-scale \\ \hline
    Koopman operators & $U_{{\rm slow},\tau}^\varepsilon$ & ~ & $U_{{\rm slow},\varepsilon t}^\varepsilon = U_{{\rm fast},t}^\varepsilon $ \\ \hline
    Koopman eigenvalues & $\lambda$ & ~ & $\varepsilon \lambda$ \\ \hline
    Koopman eigenfunctions & $\phi_\lambda^\varepsilon $ & ~ & $\phi_\lambda^\varepsilon$ \\ \hline
    \end{tabular}
    \end{center}
\end{table}

% no topic sentence
Next, we order the Koopman principal eigenvalues of the SP system (\ref{eq:vanderPol_slow}) with the periodic orbit. 
For this, the authors of Refs\cite{mishchenko2013differential,krupa2001relaxation} considered the asymptotic expansion of the periodic orbit in terms of $\varepsilon$, and we have the zeroth-order estimation of the fundamental period $T^\varepsilon = 2\pi/\omega^\varepsilon$, the associated angular frequency $\omega^\varepsilon$, and the Floquet exponent $\nu^\varepsilon$ as follows: 
\begin{equation}
    \label{eq:LCeig}
    \left.
    \begin{aligned}
        %&{\rm Period}& : & ~
        T^\varepsilon &= T^0 + o (1) \\
        %&{\rm Fundamental~angular~frequency} & : &~ 
        \omega^\varepsilon &= \omega^0 + o(1) \\
        %& {\rm Floquet ~exponent} ~& : &~ 
        \nu^\varepsilon &= \nu^0 /\varepsilon + o (\varepsilon^{-1})%,
    \end{aligned}
    \right\},
\end{equation}
where $T^0 := 3-2\ln 2$, $\omega^0 := 2\pi/T^0$, $\nu^0 := (-3/2-2\ln 2)/T^0$, and $o$ stands for the Landau notation in the sense of small $o$. 
Eq. (\ref{eq:LCeig}) shows that the two Koopman {principal} eigenvalues ${\rm i}\omega^\varepsilon$ and $\nu^\varepsilon$ of $U^\varepsilon_{{\rm slow},\tau}$ in $\tau$ scale are ordered as %by 
$\mathcal{O}(1)$ and %, 
$\mathcal{O}(\varepsilon^{-1})$, %respectively, 
where $\mathcal{O}$ is the Landau notation in the sense of big $\mathcal{O}$. 
We refer to the principal eigenvalue ${\rm i}\omega^\varepsilon$ (or $\nu^\varepsilon$) ordered $\mathcal{O}(1)$ (or $\mathcal{O}(\varepsilon^{-1})$) as a slow (or fast) eigenvalue.
Table~\ref{tab:KEV} shows the numerical results of the period $T^\varepsilon$ and slow/fast Koopman eigenvalues ${\rm i} \omega^\varepsilon$ and $\nu^\varepsilon$, verifying the accuracy of the zeroth-ordering \eref{eq:LCeig} {(see Ref\cite{use_of_Fourier_averages} and Ref\cite{phase-amplitude_reduction_nakao} for the computational methods of ${\rm i} \omega^\varepsilon$ and $\nu^\varepsilon$, respectively)}. 
\begin{table}[t]
  \caption{$\varepsilon$-dependence of period and slow/fast Koopman eigenvalues}
  \label{tab:KEV}
    \begin{center}
    \begin{tabular}{cccc}\hline
    $\varepsilon$ & $T^\varepsilon$ & ${\rm i}\omega^\varepsilon$ & $\nu^\varepsilon$\\ \hline\hline 
    1 ~& ~6.66 ~&~ i\,0.943 ~&~ $-1.06$ \\\hline
    0.1 ~& ~2.87 ~&~ i\,2.19 ~&~ $-13.3$ \\\hline
    0.01 ~& ~1.91 ~&~ i\,3.29 ~&~ $-163$ \\\hline
    \end{tabular}
    \end{center}
\end{table}

\subsection{Slow/fast Koopman eigenfunctions}
\begin{figure*}
  \centering
  \includegraphics[width=\hsize]{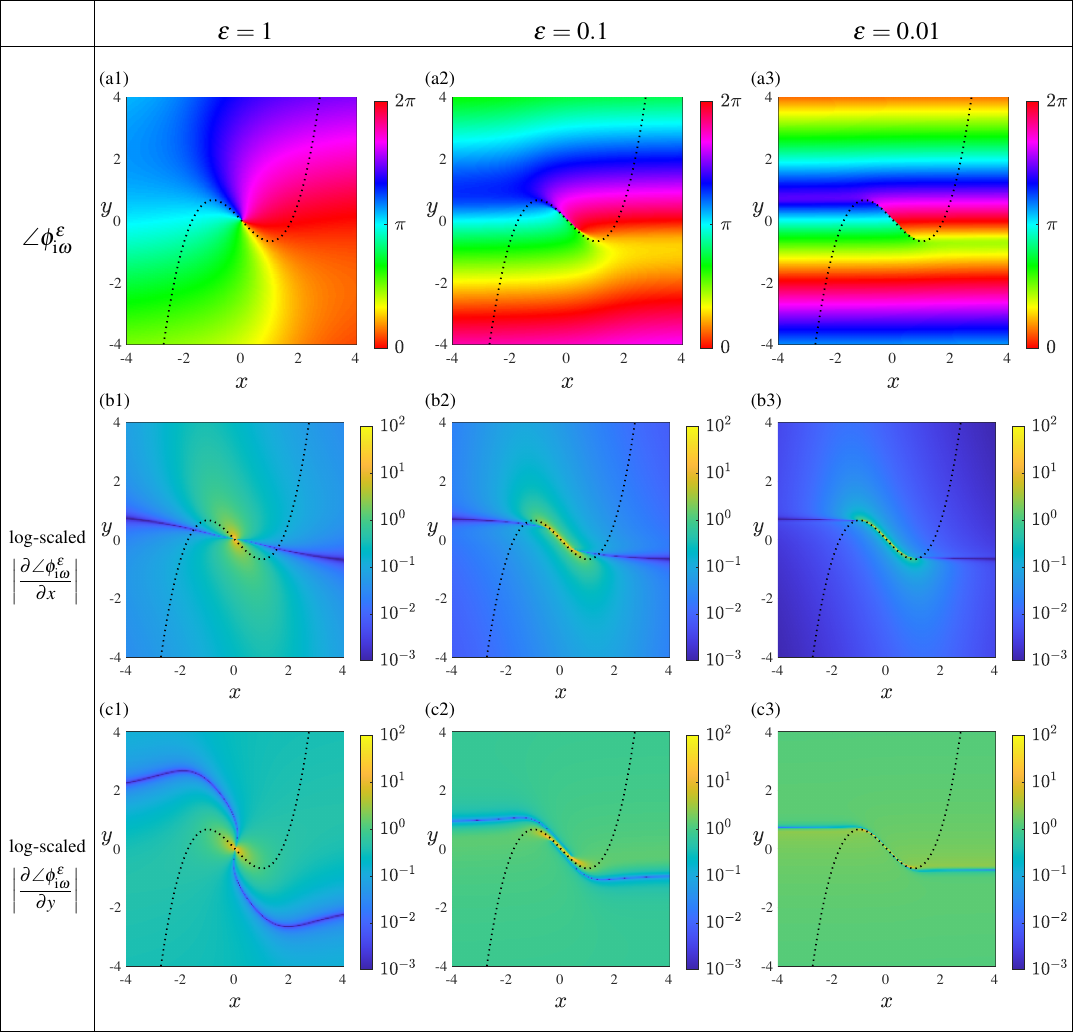}
  \caption{$\angle \phi_{{\rm i}\omega}^\varepsilon$, and log-scaled derivatives of $\phi_{{\rm i}\omega}^\varepsilon$, with $F(x,y)=0$ (dash curves).
  The derivatives are computed using finite differences.}
  \label{fig:phi_iomega}
\end{figure*}

\begin{figure*}
  \centering
  \includegraphics[width=\hsize]{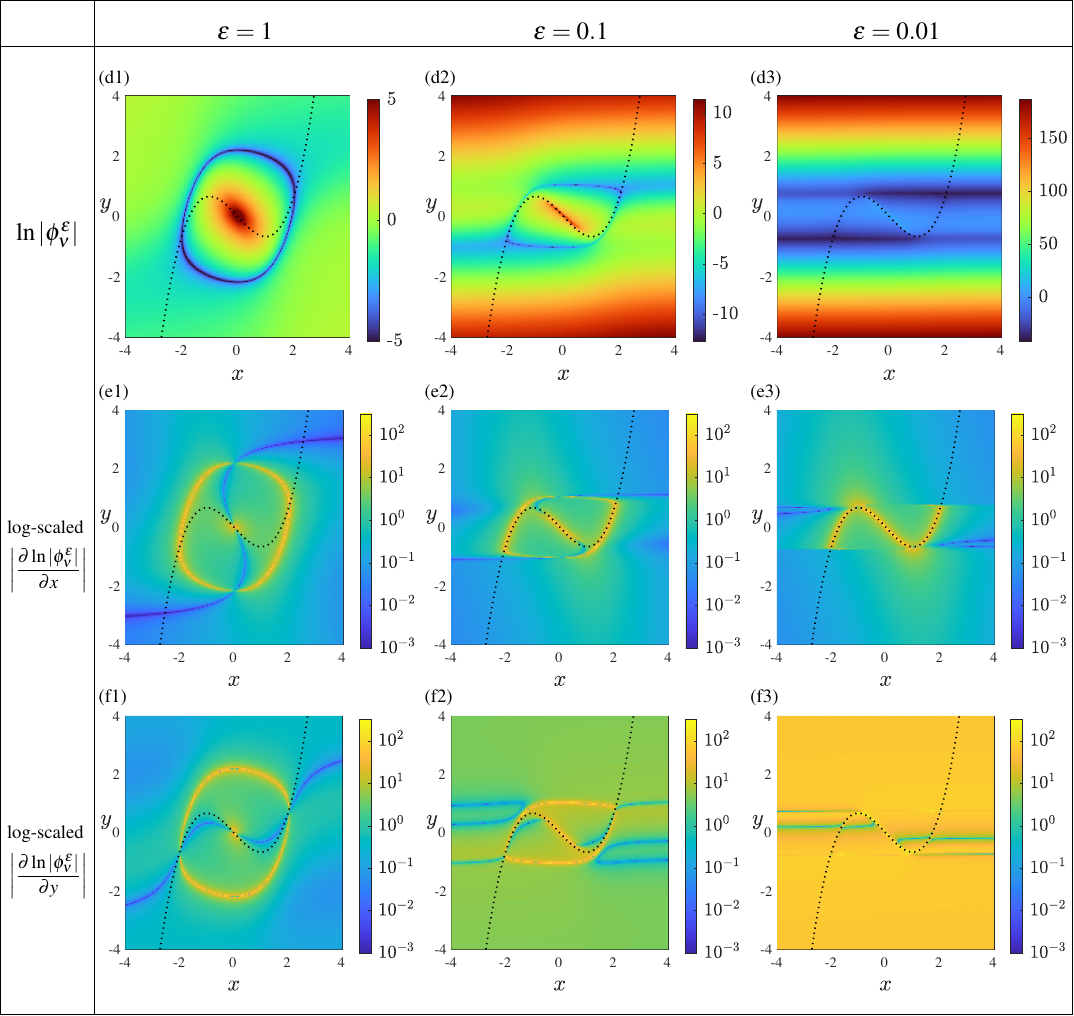}
  \caption{$\ln \phi_\nu^\varepsilon$, and log-scaled derivatives of $\phi_\nu^\varepsilon$ with $F(x,y)=0$ (dash curves).
  The derivatives are computed using finite differences.}
  \label{fig:phi_nu}
\end{figure*}

Here, we consider the Koopman {principal} eigenfunctions $\phi_{{\rm i}\omega}^\varepsilon$ and $\phi_{\nu}^\varepsilon$ associated with the eigenvalues ${\rm i}\omega^\varepsilon$ and $\nu^\varepsilon$, respectively.
We will show that as $\varepsilon$ becomes small, spectral signatures emerge as distinctive shapes of $\phi_{{\rm i}\omega}^\varepsilon$ and $\phi_{\nu}^\varepsilon$.
For this, the eigenfunctions are numerically computed for $\varepsilon = 1,~0.1$, and 0.01 with the Fourier average \cite{use_of_Fourier_averages} (or the Laplace average \cite{phase-amplitude_reduction_nakao}).
Figs.~\ref{fig:phi_iomega} and \ref{fig:phi_nu} show the numerical results including their derivatives w.r.t. the states $x$ and $y$.

\begin{figure*}
  \centering
  \includegraphics[width=\hsize]{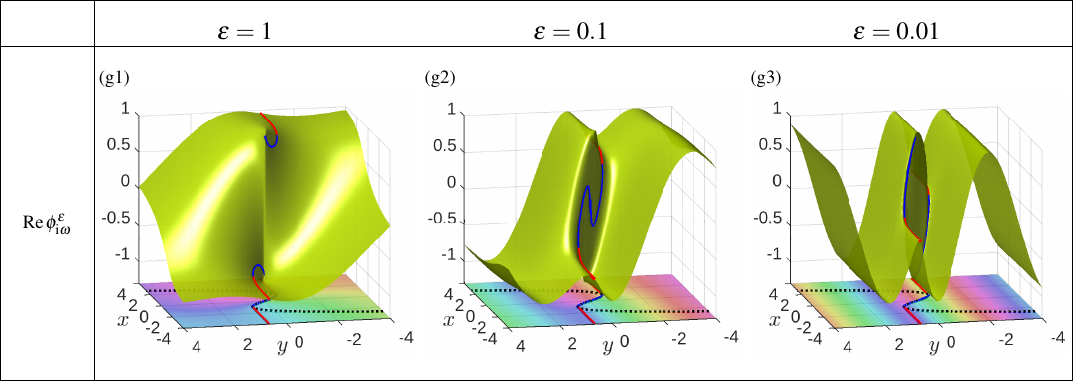}
  \caption{3D shapes of ${\rm Re} \, \phi_{{\rm i}\omega}^\varepsilon$. The 2D plots at the bottom are the same as in (a1)-(a3) of \fref{fig:phi_iomega}. In both the 3D and 2D plots, $W_0$ (or $B_\mp$) are represented by blue (or red) curves.}
  \label{fig:phi_iomega_3d}
\end{figure*}

For $\phi_{{\rm i}\omega}^\varepsilon (x,y)$ in \fref{fig:phi_iomega}, recalling the fact that $|\phi_{{\rm i}\omega}^\varepsilon|$ is constant, we depict the angles of $\phi_{{\rm i}\omega}^\varepsilon$, i.e., $\angle \phi_{{\rm i}\omega}^\varepsilon$ in (a1)-(a3) of \fref{fig:phi_iomega}.
Additionally, the (log-scaled) absolute values of its derivatives, $|\partial \angle \phi_{{\rm i}\omega}^\varepsilon /\partial x|$ and $|\partial \angle \phi_{{\rm i}\omega}^\varepsilon /\partial y|$ in (b1)-(b3) and (c1)-(c3) of \fref{fig:phi_iomega} are depicted.
For example, in (a1)-(a3), the level sets of $\angle \phi_{{\rm i}\omega}^\varepsilon$ are colored by the same color.
At first, from (a2) and (a3), it is distinctively visualized that the level sets of $\phi_{{\rm i}\omega}^\varepsilon$ are almost parallel to the $x$-axis except in the neighborhood of $W_0$.
To consider the shape of the level sets of $\phi_{{\rm i}\omega}^\varepsilon$, its derivatives provide some geometric clues. 
Especially, how the level sets are parallel to the $x$-axis is estimated by the ratio $|\partial \angle \phi_{{\rm i}\omega}^\varepsilon /\partial x | / |\partial \angle \phi_{{\rm i}\omega}^\varepsilon /\partial y|$ at each point $(x,y)$ except in the neighborhood of $W_0$, by virtue of the implicit function theorem.
The numerical results (b1)-(b3) indicate that the order of $|\partial \angle \phi_{{\rm i}\omega}^\varepsilon /\partial x |$ decreases when $\varepsilon$ changes from 1, through 0.1, to 0.01.
On the other hand, from (c1)-(c3), the order of {$|\partial \angle \phi_{{\rm i}\omega}^\varepsilon /\partial y|$} does not decrease. 
Thus, it is conjectured that as $\varepsilon$ becomes small, $|\partial \angle \phi_{{\rm i}\omega}^\varepsilon /\partial x|$ becomes smaller while $|\partial \phi_{{\rm i}\omega}^\varepsilon /\partial y|$ does not become smaller.
This conjecture $|\partial \angle \phi_{{\rm i}\omega}^\varepsilon /\partial x| \ll |\partial \angle \phi_{{\rm i}\omega}^\varepsilon /\partial y|$ leads to the distinctive shape that the level sets of $\phi_{{\rm i}\omega}^\varepsilon$ are almost parallel to the $x$-axis. 
Around $W_0$, the above observation does not hold in (a2) and (a3). 
From (b2), (b3), and (c2), (c3) in \fref{fig:phi_iomega} one can see that $|\partial \angle \phi_{{\rm i}\omega}^\varepsilon /\partial x|$ and $|\partial \angle \phi_{{\rm i}\omega}^\varepsilon /\partial y|$ increase around $W_0$ when $\varepsilon$ becomes smaller.
To understand the behavior around $W_0$, we visualize the real part ${\rm Re}\, \phi_{{\rm i}\omega}^\varepsilon$ as 3D plots shown in \fref{fig:phi_iomega_3d} (${\rm Im}\, \phi_{{\rm i}\omega}^\varepsilon$ will also provide similar visualizations, so omitted).
From (g2) and (g3) of \fref{fig:phi_iomega_3d}, ${\rm Re}\, \phi_{{\rm i}\omega}^\varepsilon$ (i.e., $\angle \phi_{{\rm i}\omega}^\varepsilon$) varies steeply around $W_0$.
This steepness is one of the distinctive shapes obtained from the numerical results.
Similarly, from (g3) of \fref{fig:phi_iomega_3d}, we see that sharp shapes emerge around $B_\mp$, which we also consider distinctive.
This sharpness can be verified numerically by computing $|\partial^2 \phi_{{\rm i}\omega}^\varepsilon /\partial y^2|$ but we omit in this manuscript.
As $\varepsilon\rightarrow 0$, it is speculated that the level sets of $\angle \phi_{{\rm i}\omega}^\varepsilon$ become completely parallel to the $x$-axis and the steep and sharp shapes become discontinuous and nonsmooth, respectively.
The latter will be investigated in the next section. 
\begin{figure*}
  \centering
  \includegraphics[width=\hsize]{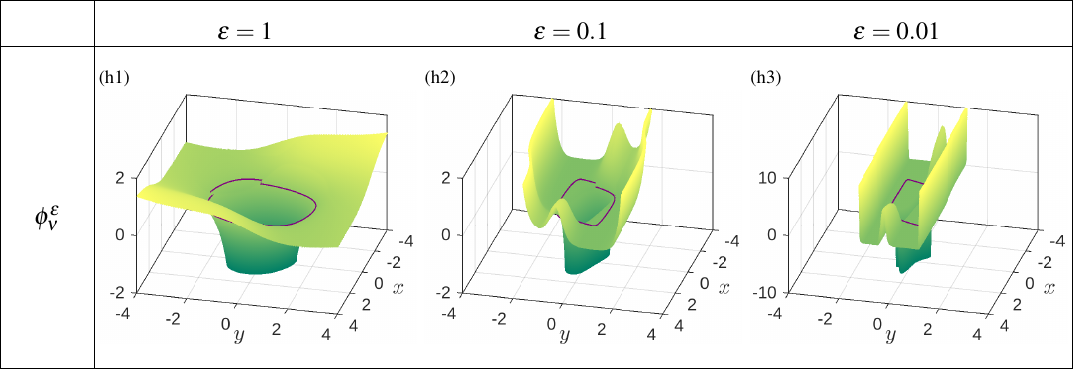}
  \caption{$\phi_\nu^\varepsilon$ with the limit cycle colored in purple.}
  \label{fig:phi_nu_3d}
\end{figure*}

For $\phi_{\nu}^\varepsilon (x,y)$ in \fref{fig:phi_nu}, because its absolute values are taken from a large range of values, 
we depict the natural logarithm of $\phi_\nu^\varepsilon$, i.e., $\ln |\phi_\nu^\varepsilon|$ in (d1)-(d3) of \fref{fig:phi_nu}.
Recalling the fact that $\phi_\nu^\varepsilon (x,y)=0$ when $(x,y)$ is on the limit cycle, one can confirm the shape of the limit cycle by tracing the the smallest part of (d1)-(d3).
%that $\ln |\phi_\nu^\varepsilon (x,y)|$ tends to $-\infty$ as $(x,y)$ approaches the limit cycle. 
Additionally, the absolute values of its derivatives, $|\partial \ln \phi_\nu^\varepsilon /\partial x|$ and $|\partial \ln \phi_\nu^\varepsilon /\partial y|$ in (e1)-(e3) and (f1)-(f3) of \fref{fig:phi_nu} are depicted. 
We also visualize $\phi_\nu^\varepsilon$ as 3D plots in \fref{fig:phi_nu_3d}.
The numerical results (d1)-(d3) exhibit that the grow rate w.r.t $y$ on the region
$\{ (x,y)\in \mathbb{R}^2~|~ |y| > 2/3\}$ becomes larger as $\varepsilon$ becomes smaller (see also (h1)-(h3) of \fref{fig:phi_nu_3d}). 
In particular, when $\varepsilon$ changes from $0.1$ to $0.01$, the grow rate increases by approximately 10th orders of magnitude, corresponding to ${\rm e}^{\nu^\varepsilon}$.
On the other hand, from (e1)-(e3) the order of $|\partial \ln |\phi_\nu^\varepsilon| /\partial x|$ does not increase when $\varepsilon$ changes from 1 to 0.1 and from 0.1 to 0.01, while from (f1)-(f3) the order of $|\partial \ln |\phi_\nu^\varepsilon| /\partial y|$ increases.
The numerical finding $|\partial \ln |\phi_\nu^\varepsilon| /\partial x| \ll |\partial \ln |\phi_\nu^\varepsilon| /\partial y|$ implies that their rapid growth emerges mostly in the $y$-direction. 
According to the implicit function theorem, the ordering of the derivatives also suggests that the level sets of $\phi_\nu^\varepsilon$ are almost parallel to the $x$-axis.
As $\varepsilon\rightarrow 0$, it is speculated that the fast Koopman eigenvalue $\nu^\varepsilon$ diverges and consequently the associated eigenfunction is not well-defined. 
The connection between the geometric properties of the SP system and the feature of $\phi_\nu^\varepsilon$ remains to be solved, which is in our future work.

\section{Singular Limit of a Koopman Eigenfunction}
In this section we aim to derive the singular limit of the Koopman eigenfunction $\phi_{{\rm i} \omega}^\varepsilon$ associated with the slow eigenvalue ${\rm i} \omega^\varepsilon$. 
Recalling that the order of the slow eigenvalues are $\mathcal{O}(1)$ from \tref{tab:KEV}, we denote the limit of ${\rm i} \omega^\varepsilon$ by ${\rm i}\omega^0$. 
The numerical result in \fref{fig:phi_iomega} suggests that there exists the singular limit of $\phi_{{\rm i} \omega}^\varepsilon$, denoted as $\phi_{{\rm i} \omega}^0$, which has discontinuous and nonsmooth sections along $W_0$ and $B_\mp$. 
Note that even though $\phi_{{\rm i} \omega}^\varepsilon \in C^\omega$ holds (see Section III), the singular limit $\phi_{{\rm i} \omega}^0$ becomes discontinuous from the numerical result.
In the following, we remove the analyticity assumption from the three functional spaces $\mathcal{F}_y$, $\bar{\mathcal{F}}$, and $\tilde{\mathcal{F}}$ which will be introduced below.
The following argument is formal, and so we do not consider any rigor formulation of the limit in the sense of normed space of observables.

To analytically derive $\phi_{{\rm i} \omega}^0$, for the fast (or slow) subsystem \eref{eq:fast_subsystem} (or \eref{eq:slow_subsystem}), we first define the observable space $\mathcal{F}_y$ (or $\bar{\mathcal{F}}$) and the Koopman operator $U_{{\rm fast},t}^0$ (or $U_{{\rm slow},\tau}^0$). 
Then, the Koopman operator $U_{{\rm fast},t}^0$ (or $U_{{\rm slow},\tau}^0$) will be defined by the flow $\boldsymbol{S}_{{\rm fast},t}^0$ (or $\boldsymbol{S}_{{\rm slow},\tau}^0$) that have been introduced in Section II. 
We also introduce the observable space $\tilde{\mathcal{F}}$ and the Koopman operator $\tilde{U}_{{\rm slow},\tau}^0$ as the composition of the flow of the singular limit (already introduced in Section II, denoted as $\tilde{\boldsymbol{S}}_{{\rm slow},\tau}^0$).
\fref{fig:observable_diagram} illustrates the relation among the five operators $U_{{\rm slow},\tau}^\varepsilon,~U_{{\rm fast},t}^\varepsilon,~U_{{\rm slow},\tau}^0,~U_{{\rm fast},t}^0$, and  $\tilde{U}_{{\rm slow},\tau}^0$,
corresponding to \fref{fig:flow_diagram}.

\begin{figure*}
    \centering
    \begin{minipage}{1\textwidth}
        {\large\[
    \xymatrix@C=0pt@R=9pt{
    { \mathcal{F}(\mathbb{R}^2)\ni f,\quad U_{{\rm slow},\tau}^\varepsilon f(x,y)}
    \ar[rr]^*{t = \tau/\varepsilon} \ar@<6ex>[dd]_*{\varepsilon \to 0} & & \ar[ll] \ar@<-8ex>[dd]_*{\varepsilon \to 0}  
    {U_{{\rm fast},t}^\varepsilon f(x,y) , \quad f \in \mathcal{F}(\mathbb{R}^2) \quad~~ }
    \\ & & \\ 
    {\bar{\mathcal{F}}(W_\pm)\ni \bar{f},\quad U_{{\rm slow},\tau}^0 \bar{f}(\bar{x},\bar{y})}
    \ar@{.}[rr] & \ar@{.>}[d]&\ar@{.}[ll] 
    {U_{{\rm fast},t}^0 f (x;y), ~~ f(\,\cdot\, ; y) \in \mathcal{F}_y (\mathbb{R})}
    \\
    & {\tilde{U}_{{\rm slow},\tau}^0 f (x,y)} &
   }
        \]}
        \caption{Five Koopman operators derived from the fast and slow subsystems, and their limits which mirror the flows' relation in \fref{fig:flow_diagram}.}
        \label{fig:observable_diagram}
    \end{minipage}
\end{figure*}

\subsection{Spectral properties of the fast subsystem}
First, we introduce the observable space and the Koopman operator for the fast subsystem, and we discuss its spectral properties.
Recalling that the fast subsystem (\ref{eq:fast_subsystem}) is a 1-dimensional system with constant $y$, i.e., a boundary layer system, the Koopman operator can be defined as a composition of the fast flow $\boldsymbol{S}_{{\rm fast},t}^0$ with an observable, both of which are 1-dimensional.  
An observable with constant $y$ is denoted as ${f}(\,\cdot\, ;y)$. 
For the discussion below, we allow observables to be discontinuous on $W_0$: we introduce $f(\,\cdot\,;y)$ and its function space $\mathcal{F}_y$ as satisfying
$$
\left. \begin{alignedat}{3}
   f(\,\cdot\,;y) \in  & C^0(\mathbb{R}) ,&\quad& |y| > 2/3 \\
   f(\,\cdot\,;y) \in  & C^0(I_-) ,&\quad& |y| \leq 2/3, ~ x\in I_- \\
   f(\,\cdot\,;y) \in  & C^0(I_+) ,&\quad& |y| \leq 2/3,~x \in I_+
\end{alignedat}
\right\} ,
$$
where $I_- = (-\infty,\gamma_0(y))$, and $I_+ = (\gamma_0(y),+\infty)$. 
Here, because the discontinuity on $(\gamma_0 (y), y)\in W_0$ is allowed, we restricted the domain of $f(\,\cdot\,;y)$ to $\mathbb{R}\backslash \{\gamma_0 (y) \}$ when $|y|\leq 2/3$.
Then, the Koopman operator of the fast subsystem is defined as $$U_{{\rm fast},t}^0 {f} (x;y):= {f} (\boldsymbol{S}_{{\rm fast},t}^0(x;y)),\quad t\geq 0.$$ 

Given that the slow eigenvalue in the fast subsystem is ${\rm i} \varepsilon \omega^\varepsilon$ and its singular limit is 0, the associated eigenfunction ${\phi}_{{\rm i}\omega}^0 (\,\cdot\,;y)$ ($y \in \mathbb{R}$) becomes a constant function from Ref\cite{isostable_isochron_and_koopmanSpectrum_fixedpoint}. 
More precisely, ${\phi}_{{\rm i}\omega}^0(x;y)$ is defined within each basin of a stable equilibrium point, excluding the basin boundary given by $\gamma_0 (y) \in W_0$. 
Therefore, ${\phi}_{{\rm i}\omega}^0 (x;y)$ becomes a constant function where $|y|>2/3$, on the other hand, becomes constant within two regions divided by $W_0$ where $|y|\leq 2/3$:
\begin{equation}
\label{eq:fastKEF}
    {\phi}_{{\rm i}\omega}^0 (x;y) = \left\{ \begin{alignedat}{3}
    &\kappa(y) ,&\quad& |y| > 2/3, \\
    &\kappa_- (y),&\quad& |y| \leq 2/3,~ x < \gamma_0 (y), \\
    &\kappa_+ (y),&\quad& |y| \leq 2/3,~ x > \gamma_0 (y). \\
\end{alignedat}
\right.
\end{equation}
where $\kappa(y)~(|y|>2/3)$ and $\kappa_\mp(y)~(|y|\leq2/3)$ are functions solely dependent on $y$.
It should be noted that 1-dimensional systems with a parameter including the boundary-layer one have been studied in Ref\cite{gaspard1995spectral}.

\subsection{Spectral properties of the slow subsystem}
Second, we introduce the observable space $\bar{\mathcal{F}}$ and the Koopman operator $U_{{\rm slow},\tau}^0$ for the slow subsystem, and we discuss its spectral properties.
The state space of the slow subsystem is $W_\mp$, and so an observable, denoted by $\bar{f}$, is defined on $W_\mp$.
Additionally, recall that the constrained flow, denoted by $\boldsymbol{S}_{{\rm slow},\tau}^0$ (see Section II), is defined on $W_\mp$ with discontinuous jumps.
The Koopman operator is defined as 
$$U_{{\rm slow},\tau}^0 \bar{f} := \bar{f}\circ {\boldsymbol{S}}_{{\rm slow},\tau}^0,\quad \tau \geq 0.$$
For discussing spectral properties of $U_{{\rm slow},\tau}^0$, it is necessary to prove the positive invariance of the observable space, i.e., 
\begin{equation}
\label{eq:slow_invariance}
    U_{{\rm slow},\tau}^0 \bar{f} \in \bar{\mathcal{F}} , \quad
\forall \tau\geq0, ~ \forall \bar{f} \in \bar{\mathcal{F}}.
\end{equation}
Here, the observable space 
$\bar{\mathcal{F}}$ is introduced as 
$$
\begin{aligned}
 \bar{\mathcal{F}} = \{ ~\bar{f}|_{W_-} \in C^0 (W_-), ~ \bar{f}|_{W_+} \in C^0 (W_+) ~|~ 
   & \lim_{ \substack{(\bar{x},\bar{y})\rightarrow (-1,2/3)\\ (\bar{x},\bar{y})\in W_-} }  \bar{f}(\bar{x},\bar{y})=\bar{f}(2, 2/3), \\
    &    \lim_{ \substack{(\bar{x},\bar{y})\rightarrow (1,-2/3)\\ (\bar{x},\bar{y})\in W_+ }}  \bar{f}(\bar{x},\bar{y})=\bar{f}(-2, -2/3) \},
\end{aligned}
$$
where $(\mp 1,\pm2/3)$ are the singular points which are closure points of $W_\pm$, and $(\pm 2, \pm 2/3) \in W_\pm$ are the drop points.
Here, We impose the continuous property $C^0$ on the observable space $\bar{\mathcal{F}}$ to show a nice spectral property of $U_{{\rm slow},\tau}^0$.
Then, the following lemma holds (see Appendices B and C for its proof).
\begin{lemma}
    $\bar{\mathcal{F}}$ is a vector space, and satisfies the positive invariance property \eqref{eq:slow_invariance}.
\end{lemma}

The $U_{{\rm slow},\tau}^0$ has an infinite countable number of Koopman eigenvalues.
One of the eigenvalues is ${\rm i}\omega^0$, as there is the periodic orbit with discontinuous jumps whose fundamental frequency is $\omega^0$. 
The eigenfunction $\bar{\phi}_{{\rm i}\omega}^0 \in \bar{\mathcal{F}}$ associated with ${\rm i}\omega^0$ is
\begin{equation}
\label{eq:slowKEFdef}
    \bar{\phi}_{{\rm i}\omega}^0 (\bar{x},\bar{y}) = \left\{ \begin{aligned}
    {\rm e}^{{\rm i}\omega^0 \varphi (\bar{x})}, &\quad \bar{x} > 1,\\
    -{\rm e}^{{\rm i}\omega^0 \varphi (|\bar{x}|)}, &\quad \bar{x} < -1,
  \end{aligned}\right.
\end{equation}
where $\varphi (\bar{x})= \ln (\bar{x}) - \bar{x}^2/2~(\bar{x}>1)$. 
This can be clarified because 
\begin{equation}
\label{eq:slowKEF}
    {U}_\tau^0 \bar{\phi}_{{\rm i}\omega}^0 (\bar{x},\bar{y})= {\rm e}^{{\rm i}\omega^0 \tau} \bar{\phi}_{{\rm i}\omega}^0 (\bar{x},\bar{y}),\quad \forall (\bar{x},\bar{y}) \in W_\mp,~\forall \tau \geq 0,
\end{equation}
holds. 
It is obvious that $\bar{\phi}_{{\rm i}\omega}^0$ is smooth on the manifolds $W_-$ and $W_+$. 
Then, we have the following lemma (see Appendix D for its proof).
\begin{lemma}
Let $\sigma (U_{{\rm slow},\tau}^0)$ denote the point spectrum (that is, the set of all Koopman eigenvalues). 
Then,
\begin{equation}
    \label{eq:eiglemma}
    \sigma (U_{{\rm slow},\tau}^0) = \{{\rm i} n \omega^0 \}_{n\in \mathbb{Z}}
\end{equation}
is satisfied. 
In addition, the associated Koopman eigenfunctions are smooth on the manifolds $W_-$ and $W_+$. 
\end{lemma}
Since the eigenfunctions of $U_{{\rm slow},\tau}^0$ are smooth on the manifolds $W_-$ and $W_+$, it is natural to equip the $r$-times {differentiable} properties ($r\geq 1$) to the observable space $\bar{\mathcal{F}}$. 
Here, we leave this for the future work because it is not straightforward to prove the positive invariance condition with $r\geq 1$ for $U_{{\rm slow},\tau}^0$.
It should be noted that
the authors of Ref\cite{shirasaka2017phase} study the isochron of the limit cycle with discontinuous jumps, which corresponds the level set of $\bar{\phi}_{{\rm i}\omega}^0$.

\subsection{Concatenation of Koopman eigenfunctions in the fast/slow subsystems}
Third, we introduce the observable space $\tilde{\mathcal{F}}$ and the Koopman operator $\tilde{U}_{{\rm slow},\tau}^0$ as the singular limit of the Koopman operator, 
thereby clarifying its spectral property.
In Section II, we introduced the state space as $\mathbb{R}^2\backslash W_0$ and defined the singular limit of the flow $\tilde{\boldsymbol{S}}_{{\rm slow},\tau}^0$ as \eref{eq:singularflow}.
An observable $\tilde{f} :\mathbb{R}^2 \backslash W_0 \to \mathbb{C}$ and its space $\tilde{\mathcal{F}}$ are defined as satisfying  
$$\left.
\begin{aligned}
    \tilde{f}|_{y = {\rm const.}} &\in {\mathcal{F}_y} \\
    \tilde{f}|_{W_\mp} &\in \bar{\mathcal{F}}  
\end{aligned}\right\}.
$$
In other words, the restriction of $\tilde{f}$ on $y={\rm const.}$ is introduced as an observable of the fast subsystem, and the restriction of $\tilde{f}$ on $W_\mp$ is introduced as an observable of the slow subsystem.
Then, the Koopman operator is defined as 
$$\tilde{U}_{{\rm slow},\tau}^0 \tilde{f} = \tilde{f}\circ \tilde{\boldsymbol{S}}_{{\rm slow},\tau}^0,\quad \tau \geq 0,$$
and it can be seen that
$$
\tilde{U}_{{\rm slow},\tau}^0 \tilde{f} \in \tilde{\mathcal{F}}, \quad
\forall \tau\geq0, ~ \forall \tilde{f} \in \tilde{\mathcal{F}}
$$ 
holds.
The Koopman eigenvalues and eigenfunctions of $\tilde{U}_{{\rm slow},\tau}^0$ are defined as 
$$
\tilde{U}_{{\rm slow},\tau}^0 \tilde{\phi}_{\lambda} (x,y)= {\rm e}^{\lambda \tau} \tilde{\phi}_{\lambda} (x,y),~\forall (x,y)\in \mathbb{R}^2\backslash W_0,~\forall \tau\geq 0.
$$
Since the observable $\tilde{f}$ is introduced as the concatenation of {those of} the fast and slow subsystems, the Koopman eigenfunction $\tilde{\phi}_{{\rm i}\omega}^0$ associated with ${\rm i}\omega^0$ is similarly derived:
$$
\tilde{\phi}_{{\rm i} \omega}^0 (x,y) = \left\{ 
  \begin{aligned}
    {\rm e}^{{\rm i}\omega^0 \varphi (\gamma_+(y))}&,\quad (x,y)\in D_+, \\ 
    -{\rm e}^{{\rm i}\omega^0 \varphi (|\gamma_-(y)|)}&,\quad (x,y)\in D_-.
  \end{aligned}
\right.
$$
One can see that
$\tilde{\phi}_{{\rm i}\omega}^0 \in \tilde{\mathcal{F}}$ holds since $ \tilde{\phi}_{{\rm i}\omega}^0|_{W_\mp}$ satisfies Eq. (\ref{eq:slowKEF}) and $\tilde{\phi}_{{\rm i}\omega}^0|_{y={\rm const.}}$ satisfies Eq. (\ref{eq:fastKEF}). 
It is obvious that $\bar{\phi}_{{\rm i}\omega}^0$ is smooth on the regions $D_-$ and $D_+$.
The following theorem provides a spectral characterization of the singular limit (its proof is straightforward from the two Lemmas, so omitted). 
\begin{theorem}
Let $\sigma (\tilde{U}_{{\rm slow},\tau}^0)$ denote the point spectrum (that is, the set of all Koopman eigenvalues). 
Then,
\begin{equation}
    \label{eq:eigtheorem}
    \sigma (\tilde{U}_{{\rm slow},\tau}^0) = \{{\rm i} n \omega^0 \}_{n\in \mathbb{Z}}
\end{equation}
is satisfied. 
In addition, the associated Koopman eigenfunctions are smooth on the regions $D_-$ and $D_+$. 
\end{theorem}

\begin{figure*}[t]
  \begin{tabular}{cc}
    \begin{minipage}[t]{0.4\hsize}
      \centering
      \includegraphics[keepaspectratio, scale=0.8]{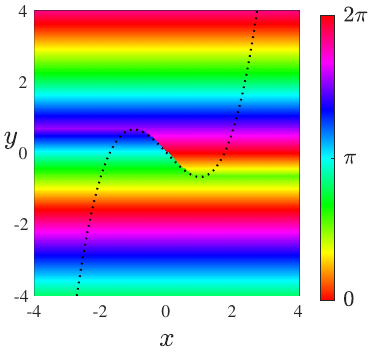}
      \subcaption{$\angle \tilde{\phi}_{{\rm i} \omega}^0$.}
    \end{minipage}&
    \begin{minipage}[t]{0.4\hsize}
      \centering
      \includegraphics[keepaspectratio, scale=0.8,page=1]{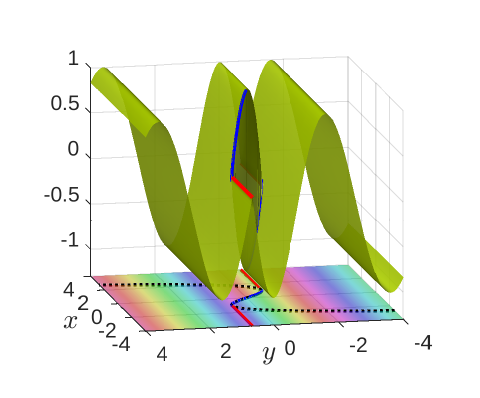}
      \subcaption{${\rm Re} \, \tilde{\phi}_{{\rm i} \omega}^0$.}
    \end{minipage}
  \end{tabular}
  \caption{Singular limit of the Koopman eigenfunction $\phi_{{\rm i} \omega}^\varepsilon$.}
  \label{fig:singular}
\end{figure*}
\fref{fig:singular} plots $\tilde{\phi}_{{\rm i}\omega}^0$.
From the geometrical perspective, the Koopman eigenfunction $\tilde{\phi}_{{\rm i}\omega}^0$ behaves in a manner similar to the distinct shapes observed in Sec.~III;
but it becomes discontinuous and nonsmooth on $W_0$ and $B_\mp$, respectively.
These capture the properties of the singular limit of the flow $\tilde{\boldsymbol{S}}_{{\rm slow},\tau}^0$ as follows: 
the level set including $p\in W_\mp$ coincide with the {\em fiber} $\pi^{-1} (p)$, or, in other words, the level sets coincide with the leaves of the foliation of the attracting critical manifold $W_\mp$ (check the precise definition in Ref \cite{takens1976constrained}).
Additionally, the discontinuous and nonsmooth sections correspond to the discontinuous section of $\pi$.
\tref{tab:flow_KEF} summarizes the correspondences between $\tilde{\phi}_{{\rm i}\omega}^0$ (or ${\phi}_{{\rm i}\omega}^\varepsilon$) and the geometric properties of $\tilde{\boldsymbol{S}}_{{\rm slow},\tau}^0$ (or {$\boldsymbol{S}_{{\rm slow},\tau}^\varepsilon$}). 
The center column of \tref{tab:flow_KEF} shows that $\tilde{\phi}_{{\rm i}\omega}^0$ captures three properties of $\tilde{\boldsymbol{S}}_{{\rm slow},\tau}^0$;
(1.) the foliation of the critical manifold, (2.) discontinuity of $\pi$ on $W_0$, and (3.) discontinuity of $\pi$ on $B_\mp$.

{Note that the authors of Ref\cite{govindarajan2016operator} develop the Koopman operator for piecewise smooth dynamical systems. 
We speculate that our result of the singular limit contributes a new development of the Koopman operator for piecewise smooth dynamical systems. }
\begin{table}
    \caption{(Center column) Correspondences between the singular limit of the flow and the singular limit of the eigenfunction. (Right column) Conjecture of the relation between the SP flow and the SP eigenfunction.
    }
    \label{tab:flow_KEF}
    \hspace{-12.5mm}
    \begin{tabular}{|c|c|c|}\hline
        & singular limit & $\varepsilon>0$  \\ \hline
       flow &  $ \tilde{\boldsymbol{S}}_{{\rm slow},\tau}^0 :\left\{\begin{aligned}
        &1.~\text{foliation of the critical manifold}\\ &2.~\text{discontinuity of $\pi$ on $W_0$} \\ &3.~\text{discontinuity of $\pi$ on $B_\mp$}
    \end{aligned}\right.$ &  $ \boldsymbol{S}_{{\rm slow},\tau}^\varepsilon :\left\{\begin{aligned}
        &1.~\text{foliation of slow manifold} \\ &2.~\text{non-existence region of %} \\ & \qquad \text{
        foliation around $W_0$} \\ &3.~\text{non-existence region of %} \\ & \qquad \text{
        foliation around $B_\mp$}
    \end{aligned}\right.$ \\ \hline
       $\begin{aligned}
           &\text{Koopman} \\ &\text{eigenfunction}
       \end{aligned}$ & $ \tilde{\phi}_{{\rm i}\omega}^0 :\left\{\begin{aligned}
        &1.~\text{completely parallel level sets} \\ &2.~\text{discontinuity on $W_0$} \\ &3.~\text{nonsmoothness on $B_\mp$}
    \end{aligned}\right.$ &  $ {\phi}_{{\rm i}\omega}^\varepsilon :\left\{\begin{aligned}
        &1.~\text{almost parallel level sets} \\ &2.~\text{steepness around $W_0$} \\ &3.~\text{sharpness around $B_\mp$}
    \end{aligned}\right.$ \\ \hline
    \end{tabular}
\end{table}

\subsection{Discussion}
Here, based on the properties of the Koopman eigenfuntion $\tilde{\phi}_{{\rm i}\omega}^0$, that is, the center column of \tref{tab:flow_KEF}, we infer the connection between the flow of the SP system (denoted by $\boldsymbol{S}_{{\rm slow},\tau}^\varepsilon$) and the Koopman eigenfunction $\phi_{{\rm i}\omega}^\varepsilon$.
In Section III we showed that $\phi_{{\rm i} \omega}^\varepsilon$ exhibits the distinct shapes: the level sets are almost parallel to the $x$-axis, and steep (or sharp) shapes emerge around $W_0$ (or $B_\mp$), see \fref{fig:phi_iomega}.
As we change $\varepsilon$ to 0, the level sets become completely parallel.
This change is similar to the change from the SP flow to the singular limit: from {\em global foliation of slow manifold} \cite{eldering2018global} to the foliation of the critical manifold (i.e., fibers $\pi^{-1}$). 
From the analogy of the connections between $\tilde{\boldsymbol{S}}_{{\rm slow},\tau}^0$ and $\tilde{\phi}_{{\rm i}\omega}^0$,
it is speculated that the level sets of $\phi_{{\rm i} \omega}^\varepsilon$ coincide with the leaves of the foliation of the slow manifold.
The authors of Ref\cite{isostable_isochron_and_koopmanSpectrum_fixedpoint} show for a system with an asymptotically stable equilibrium point, level sets of Koopman eigenfunctions correspond to fibers of a slow manifold.
It is also speculated that the steep and sharp shapes of $\phi_{{\rm i} \omega}^\varepsilon$ emerge on regions where no foliation of slow manifold exists.
For sufficiently small $\varepsilon$, the regions are around $W_0$, $B_\mp$.
The discussed link between $\boldsymbol{S}_{{\rm slow},\tau}^\varepsilon$ and $\phi_{{\rm i} \omega}^\varepsilon$ is summarized in right column of \tref{tab:flow_KEF}.
This discussion has two difficulties: 
First, the foliation of the slow manifold cannot be globally extended since the slow manifolds formed along $W_\mp$ are not {\em inflowing invariant}, see Ref\cite{eldering2018global} for its detail discussion.
Second, it has not been asserted that level sets of Koopman eigenfunctions correspond to the fibers of the slow manifold for the SP van der Pol system, due to the lack of inflowing invariance.

Finally, we generalize the Koopman analysis of the SP van der Pol oscillator to a planar SP oscillator
\begin{equation}
    \label{eq:generalSP}
    \left\{\begin{aligned}
        \varepsilon \dot{x} &= F(x,y), \\
        \dot{y} &= G(x,y),
    \end{aligned}\right. 
\end{equation}
where $F,~G:\mathbb{R}^2 \to \mathbb{R}$ are analytic functions and $\varepsilon$ is a positive small parameter.
We assume that a relaxation oscillation with period $T^\varepsilon = T^0 + o(1)$ emerges on a basin $\mathcal{B}^\varepsilon$ (see Ref\cite{mishchenko2013differential} for detailed assumptions regarding the existence of a relaxation oscillation).
One can consider the flow $\boldsymbol{S}_{{\rm slow},\tau}^\varepsilon$ with $0< \varepsilon \ll 1$ and also consider its singular limit $\tilde{\boldsymbol{S}}_{{\rm slow},\tau}^0$ on the basin $\mathcal{B}^0$ by introducing the projection $\pi$ to the critical manifold and the constrained flow. 
From the analogy of our discussion for the SP van der Pol system (Section III), the Koopman operators $U_{{\rm slow},\tau}^\varepsilon$ is defined as 
$$
U_{{\rm slow},\tau}^\varepsilon f := f \circ \boldsymbol{S}_{{\rm slow},\tau}^\varepsilon,\quad {\forall \tau \geq 0}, 
$$ 
where $f:\mathcal{B}^\varepsilon \to \mathbb{C}$, $f\in \mathcal{F}$ is an observable and $\mathcal{F}$ is a set of analytic functions. 
Then, there exists a Koopman eigenfunction $\phi_{{\rm i}\omega}^\varepsilon \in \mathcal{F}$ satisfying
$$
U_{{\rm slow},\tau}^\varepsilon \phi_{{\rm i}\omega}^\varepsilon (x,y) = {\rm e}^{{\rm i}\omega^\varepsilon \tau} \phi_{{\rm i}\omega}^\varepsilon(x,y),\quad \forall (x,y) \in \mathcal{B}^\varepsilon ,~{\forall \tau \geq 0},
$$
where $\omega^\varepsilon = 2\pi / T^\varepsilon$ is the angular frequency of the relaxation oscillation.
Also, from the analogy of our discussion for the singular limit of the SP van der Pol system (Section IV), one can define the singular limit of the Koopman operator as
$$
\tilde{U}_{{\rm slow},\tau}^0 \tilde{f} := \tilde{f} \circ \tilde{\boldsymbol{S}}_{{\rm slow},\tau}^0,\quad {\forall \tau \geq 0},
$$
where $\tilde{f}:\mathcal{B}^0 \to \mathbb{C}$, $\tilde{f}\in \tilde{\mathcal{F}}$ is an observable and $\tilde{\mathcal{F}}$ a set of piecewise continuous functions. 
Then, there exists a Koopman eigenfunction $\tilde{\phi}_{{\rm i}\omega}^0 \in \mathcal{F}$ satisfying
$$
\tilde{U}_{{\rm slow},\tau}^0 \tilde{\phi}_{{\rm i}\omega}^0 (x,y) = {\rm e}^{{\rm i}\omega^0 \tau} \tilde{\phi}_{{\rm i}\omega}^0 (x,y),\quad \forall (x,y) \in \mathcal{B}^0 ,~ {\forall \tau \geq 0},
$$
where $\omega^0 = 2\pi / T^0$ is the zeroth approximation of the angular frequency $\omega^\varepsilon$. 
From the geometrical perspective, 
$\tilde{\phi}_{{\rm i}\omega}^0$ captures the properties of $\tilde{\boldsymbol{S}}_{{\rm slow},\tau}^0$ as follows: the level sets of $\tilde{\phi}_{{\rm i}\omega}^0$ (including $\pi$) coincide with the leaves of the foliation of the attracting critical  manifold (specifically, the leaves of the foliation of $W_\mp$ in the SP van der Pol system). 
Additionally, discontinuous and nonsmooth sections of $\tilde{\phi}_{{\rm i}\omega}^0$ can emerge on the sections where $\pi$ becomes discontinuous.
By changing from $\varepsilon =0$ to $0<\varepsilon \ll 1$, it can be conjectured that $\phi_{{\rm i}\omega}^\varepsilon$ has distinct shapes mirroring geometric properties of $\boldsymbol{S}_{{\rm slow},\tau}^\varepsilon$ as follows:
the level sets of $\phi_{{\rm i}\omega}^\varepsilon$ coincide with the leaves of the foliation of the attracting slow manifold. 
Additionally, steep and sharp sections of ${\phi}_{{\rm i}\omega}^\varepsilon$ emerge on the sections where no foliation of the attracting slow manifold exists.
{Extending our approach to high-dimensional SP systems with complex dynamics is in our future work.}

%%%%%%%%%%%%%%%%%%%%%%%%%%%%%%%%%%%%%%%%%%%%%%%%%%%%%%%%%%%
\section{Conclusion}
This paper performed the Koopman analysis of the SP van der Pol system. 
We showed that as $\varepsilon$ smaller, the spectral signatures emerge in the distinctive shapes of the Koopman {principal} eigenfunctions. 
We also analyzed the Koopman operator for the slow and fast subsystems, and its singular limit, and their spectral properties.
From the analysis, it is inferred that the distinctive shapes mirror the geometric properties of the SP flow.
We believe that the Koopman operator for SP systems has the potential to construct a novel perturbation theory as asymptotic expansions of the operator itself, observables evolution, and flow. 
One of the future works is 
%theoretical guarantees for this inference.
the rigorous arguement with normed spaces for the Koopman operators of $\varepsilon>0$ and $\varepsilon=0$. 

\section*{Acknowledgment}
We appreciate to Dr.~Alexandre Mauroy for his careful reading and valuable comments of the manuscript. 
The work was partially supported by JSPS KAKENHI Grant Number 23H01434 and JST Moonshot R\&D Grant Number JPMJMS2284.

\appendix

\section{The constrained flow}
We construct the constrained flow of the semi-explicit DAE (\ref{eq:slow_subsystem}) induced by the SP van der Pol oscillator. 
The solution of an initial state $(\bar{x}_0,\bar{y}_0)\in W_+$ of the semi-explicit DAE (\ref{eq:slow_subsystem}) in $\tau\in[0,\tau_{\rm s})$ is constructed as
\begin{equation}
  \label{eq:sol}
  \bar{x}(\tau) = \varphi^{-1} (\varphi (\bar{x}_0) + \tau), \quad \bar{y}(\tau) = \bar{x}(\tau)^3/3 - \bar{x}(\tau),
\end{equation}
where $\varphi (\bar{x}) = \ln |\bar{x}| - {\bar{x}^2}/{2}$ with $\bar{x}>1$ and $\tau_{\rm s}$ is the time {when the solutions reach one of the singular points ${\rm J}_+$ or ${\rm J}_-$ for the first time, satisfying $\tau_{\rm s} = \varphi (1) - \varphi(\bar{x}_0)$}.
At $\tau=\tau_{\rm s}$,  ${\rm J}_+$, the solution jumps discontinuously to the drop point $(-2,-2/3)$, so the solution extended to $\tau\in [0,\tau_{\rm s}+T^0/2)$ is obtained as 
\begin{equation}
  \nonumber
  \begin{aligned}
      \bar{x}(\tau) &= \left\{ 
    \begin{aligned}
      \varphi^{-1} (\varphi (\bar{x}_0) + \tau) & \quad \tau\in [0,\tau_{\rm s}), \\
      - \varphi^{-1} (\varphi (|-2|) + \tau) & \quad \tau\in [\tau_{\rm s},\tau_{\rm s}+T^0/2),
    \end{aligned} 
  \right. \\
   y(\tau) &= x(\tau)^3/3 - x(\tau), 
  \end{aligned} 
\end{equation}
where $T^0/2 = 3/2 - \ln 2$ is the time it takes to move from the drop point to the singular point. 
Repeating the time extension iteratively, the solution in $\tau \in [0,\infty)$ is obtained. 
The solution of an initial state $(\bar{x}_0,\bar{y}_0) \in W_-$ is also constructed but omitted. 
As a result, the solution {$(\bar{x}_\tau(\bar{x}_0,\bar{y}_0),\bar{y}_\tau(\bar{x}_0,\bar{y}_0))$} of the semi-explicit DAE (\ref{eq:slow_subsystem}) is constructed as follows:
\begin{equation}
  \label{eq:sol3}
  \begin{aligned}
    \bar{x}_\tau (\bar{x}_0,\bar{y}_0) &= \left\{ 
    \begin{aligned}
      &{\rm sgn} (\bar{x}_0) \, \varphi^{-1} (\varphi (\bar{x}_0) + \tau),  \quad \tau\in [0,\tau_{\rm s}), \\
     &- {\rm sgn} (\bar{x}_0) \, \varphi^{-1} (\varphi (2) + \tau),  \quad \tau \in \bigcup_{n=0}^{\infty}[\tau_{\rm s} + nT^0,\tau_{\rm s}+(n+1/2)T^0), \\
     &{\rm sgn} (\bar{x}_0) \, \varphi^{-1} (\varphi (2) + \tau),  \quad \tau \in \bigcup_{n=0}^{\infty}[\tau_{\rm s} + (n+1/2)T^0,\tau_{\rm s}+(n+1)T^0),
    \end{aligned}
    \right. \\
    \bar{y}_\tau (\bar{x}_0,\bar{y}_0) &= \bar{x}_\tau (\bar{x}_0,\bar{y}_0)^3/3 - \bar{x}_\tau (\bar{x}_0,\bar{y}_0).
  \end{aligned}
\end{equation}
Thus, the singular limit of the flow of (\ref{eq:slow_subsystem}) constrained by $W_\pm$ is defined as 
$${\boldsymbol{S}}_{{\rm slow},\tau}^0 (\bar{x},\bar{y}) := (\bar{x}_\tau (\bar{x},\bar{y}),\bar{y}_\tau (\bar{x},\bar{y})),\quad \forall(\bar{x},\bar{y})\in W_\pm.$$ 
Obviously, $\{{\boldsymbol{S}}_{{\rm slow},\tau}^0 : W_+ \cup W_- \to W_+ \cup W_- ,~\tau \geq 0\}$ is a semi-group, i.e.,
$$\left. \begin{aligned}
    {\boldsymbol{S}}_{{\rm slow},\tau_1 + \tau_2}^0 (\bar{x},\bar{y}) &= {\boldsymbol{S}}_{{\rm slow},\tau_1 }^0 \circ {\boldsymbol{S}}_{{\rm slow},\tau_1 }^0 (\bar{x},\bar{y})  \\
     {\boldsymbol{S}}_{{\rm slow},0}^0 (\bar{x},\bar{y}) &= (\bar{x},\bar{y})
\end{aligned}
\right\} \quad \forall (\bar{x},\bar{y})\in W_\mp .$$

\section{Proof that $\bar{\mathcal{F}}$ is a vector space}
\begin{proof}
Let us define the addition of observables as $(\alpha \bar{f} + \beta \bar{g})(\bar{x},\bar{y}) := \alpha \bar{f}(\bar{x},\bar{y}) + \beta \bar{g}(\bar{x},\bar{y}) $ for $\bar{f},\bar{g}\in \bar{\mathcal{F}}$ and $\alpha,\beta\in\mathbb{C}$. 
The zero vector is identically the zero-valued function on $W_\pm$, and an inverse element of $\bar{f}$ is $(-\bar{f})(\bar{x},\bar{y}):= -\bar{f}(\bar{x},\bar{y})$.
Obviously, the inclusion property 
$ \alpha \bar{f} + \beta \bar{g} \in C^0 (W_\mp)$ holds. 
At the singular point $(-1,2/3)$, we estimate the limit
$$\begin{aligned}
    &\lim_{\substack{(\bar{x},\bar{y})\to (-1,2/3)\\ (\bar{x},\bar{y})\in W_+}} (\alpha \bar{f} + \beta \bar{g}) (\bar{x},\bar{y}) \\&~~= \alpha \lim_{\substack{(\bar{x},\bar{y})\to (-1,2/3)\\ (\bar{x},\bar{y})\in W_+}}  \bar{f} (\bar{x},\bar{y}) + \beta \lim_{\substack{(\bar{x},\bar{y})\to (-1,2/3)\\ (\bar{x},\bar{y})\in W_+}} \bar{g} (\bar{x},\bar{y}) \\&~~
        = \alpha \bar{f}(2,2/3) + \beta \bar{g}(2,3/2) = (\alpha \bar{f} + \beta \bar{g})(2,2/3).
\end{aligned}
$$
The limit also holds for the singular point $(+1,-2/3)$. 
This shows that $\bar{\mathcal{F}}$ is closed under the addition and scalar multiplication.
\end{proof}

\section{Proof of \eref{eq:slow_invariance}}
\begin{proof}
We will prove that for all $\tau \geq 0$ and $\bar{f}\in \bar{\mathcal{F}}$, 
\begin{equation}
\label{eq:UC0}
    U_{{\rm slow},\tau}^0 \bar{f} \in C^0 (W_\mp),
\end{equation}
\begin{equation}
\label{eq:Usingularp}
    \lim_{ \substack{(\bar{x},\bar{y})\rightarrow (-1,2/3)\\ (\bar{x},\bar{y})\in W_-} } U_{{\rm slow},\tau}^0 \bar{f}(\bar{x},\bar{y})= U_{{\rm slow},\tau}^0 \bar{f}(2, 2/3),
\end{equation}
\begin{equation}
\label{eq:Usingularm}
    \lim_{ \substack{(\bar{x},\bar{y})\rightarrow (1,-2/3)\\ (\bar{x},\bar{y})\in W_+} } U_{{\rm slow},\tau}^0 \bar{f}(\bar{x},\bar{y})= U_{{\rm slow},\tau}^0 \bar{f}(-2, -2/3),
\end{equation}
using the properties of $\bar{\mathcal{F}}$:
\begin{equation}
\label{eq:C0}
    \bar{f} \in C^0 (W_\mp),
\end{equation}
\begin{equation}
\label{eq:singularp}
    \lim_{ \substack{(\bar{x},\bar{y})\rightarrow (-1,2/3)\\ (\bar{x},\bar{y})\in W_-} } \bar{f}(\bar{x},\bar{y})= \bar{f}(2, 2/3),
\end{equation}
\begin{equation}
\label{eq:singularm}
    \lim_{ \substack{(\bar{x},\bar{y})\rightarrow (1,-2/3)\\ (\bar{x},\bar{y})\in W_+} } \bar{f}(\bar{x},\bar{y})= \bar{f}(-2, -2/3).
\end{equation}
In other words, we will derive (C4)-(C6) $\Longrightarrow$ (C1)-(C3).
If $\tau = 0$, (C1)-(C3) are equivalent to (C4)-(C6), thus we only consider the case $\tau > 0$.
To show \eqref{eq:UC0},
we show $U_{{\rm slow},\tau}^0 \bar{f} (\bar{x},\bar{y})$ is continuous at any state $(\bar{x},\bar{y})\in W_+$ ($(\bar{x},\bar{y})\in W_-$ can be shown similarly and is omitted). 
We introduce the set $\{\tau_{\rm s} + nT^0/2\}_{n\in \mathbb{N}_{\geq 0}}$ ($\mathbb{N}_{\geq 0}$: nonnegative integers) at which discontinuous jumps emerge (see Appendix A), where $\tau_{\rm s} = \varphi (1) - \varphi (\bar{x})$ and $\varphi (x) = \ln |x| - x^2/2$ with $x>1$.
Let us consider the following three cases:
$$
\begin{aligned}
    &\tau \notin \{\tau_{\rm s} + nT^0/2\}_{n\in \mathbb{N}_{\geq 0}}, \\
    &\tau \in \{\tau_{\rm s} + nT^0/2\}_{n=0,2,4,\ldots},\\
    &\tau \in \{\tau_{\rm s} + nT^0/2\}_{n=1,3,5,\ldots}.
\end{aligned}
$$
If $\tau\notin \{\tau_{\rm s} + nT^0/2\}_{n\in \mathbb{N}_{\geq 0}}$, then there exists the vicinity $V\subset W_+$ of $(\bar{x},\bar{y})$ such that $\boldsymbol{S}_{{\rm slow},\tau}^0 |_{V}$ is a continuous map. 
This, along with \eqref{eq:C0}, implies that $U_{{\rm slow},\tau}^0 \bar{f} (\bar{x},\bar{y})$ is continuous.
If $\tau \in \{\tau_{\rm s} + nT^0/2\}_{n=0,2,4,\ldots}$, then, by considering that $W_+$ is 1-dimensional (and \eqref{eq:C0}), it is possible to calculate the limit from the right part $\{(\bar{x}',\bar{y}')\in W_+~|~ \bar{x}'>\bar{x}\}$
$$
\begin{aligned}
     \lim_{ \substack{(\bar{x}',\bar{y}')\rightarrow (\bar{x},\bar{y})\\ (\bar{x}',\bar{y}')\in W_+,~ \bar{x}'>\bar{x}} } U_{{\rm slow},\tau}^0 \bar{f} (\bar{x}',\bar{y}') 
     &= \lim_{ \substack{(\bar{x}',\bar{y}')\rightarrow (\bar{x},\bar{y})\\ (\bar{x}',\bar{y}')\in W_+,~ \bar{x}'>\bar{x}} } \bar{f} (\boldsymbol{S}_{{\rm slow},\tau}^0 (\bar{x}',\bar{y}')) \\
     &= \lim_{ \substack{(\bar{x},\bar{y})\rightarrow (1,-2/3)\\ (\bar{x},\bar{y})\in W_+} } \bar{f}(\bar{x},\bar{y})\quad (\,\because~\text{\eref{eq:sol3}}\,),
\end{aligned}
$$
and the limit from the left part $\{(\bar{x}',\bar{y}')\in W_+~|~ \bar{x}'<\bar{x}\}$
$$
\begin{aligned}
     \lim_{ \substack{(\bar{x}',\bar{y}')\rightarrow (\bar{x},\bar{y})\\ (\bar{x}',\bar{y}')\in W_+,~ \bar{x}'<\bar{x}} } U_{{\rm slow},\tau}^0 \bar{f} (\bar{x}',\bar{y}') 
     &= \lim_{ \substack{(\bar{x}',\bar{y}')\rightarrow (\bar{x},\bar{y})\\ (\bar{x}',\bar{y}')\in W_+,~ \bar{x}'<\bar{x}} } \bar{f} (\boldsymbol{S}_{{\rm slow},\tau}^0 (\bar{x}',\bar{y}')) \\
     &= \bar{f}(-2,-2/3) \quad (\,\because~\text{\eref{eq:sol3}}\,).
\end{aligned}
$$
From \eqref{eq:singularm}, 
both the limits are the same, thus $U_{{\rm slow},\tau}^0 \bar{f}$ is continuous.
It can be shown similarly if $\tau \in \{\tau_{\rm s} + nT^0/2\}_{n=1,3,5,\ldots}$.
Therefore, \eqref{eq:UC0} holds for all $\tau> 0$ and $\bar{f}\in \bar{\mathcal{F}}$.

Next, we prove \eqref{eq:Usingularm}. 
At first, we have 
$$
\begin{aligned}
    \boldsymbol{S}_{{\rm slow},\tau}^0 (\bar{x},\bar{y}) &= \boldsymbol{S}_{{\rm slow},\tau-\tau_{\rm s}(\bar{x})}^0 \circ \boldsymbol{S}_{{\rm slow},\tau_{\rm s}(\bar{x})}^0 (\bar{x},\bar{y}) \\
 &= \boldsymbol{S}_{{\rm slow},\tau-\tau_{\rm s}(\bar{x})}^0 (-2,-2/3), \quad \forall (\bar{x},\bar{y})\in W_+,
\end{aligned}
$$
where $\tau_{\rm s} (\bar{x}) = \varphi(1) - \varphi(\bar{x})$. 
Note that 
\begin{equation}
\label{eq:taus0}
    \lim_{\bar{x}\to 1} \tau_{\rm s} (\bar{x}) = 0.
\end{equation}
The left-hand side of \eqref{eq:Usingularm} can be represented as 
$$
\begin{aligned}
\lim_{ \substack{(\bar{x},\bar{y})\rightarrow (1,-2/3)\\ (\bar{x},\bar{y})\in W_+} } U_{{\rm slow},\tau}^0 \bar{f}(\bar{x},\bar{y}) 
   &= \lim_{ \substack{(\bar{x},\bar{y})\rightarrow (1,-2/3)\\ (\bar{x},\bar{y})\in W_+} } \bar{f} ( \boldsymbol{S}_{{\rm slow},\tau}^0 (\bar{x},\bar{y})) \\
   &= \lim_{ \substack{(\bar{x},\bar{y})\rightarrow (1,-2/3)\\ (\bar{x},\bar{y})\in W_+} } \bar{f} ( \boldsymbol{S}_{{\rm slow},\tau-\tau_{\rm s}(\bar{x})}^0 {(-2,-2/3)})
\end{aligned}
$$
Let us consider the following three cases:
$$
\begin{aligned}
    &\tau \notin \{nT^0/2\}_{n\in \mathbb{N}_{\geq 0}}, \\
    &\tau \in \{nT^0/2\}_{n=1,3,5,\ldots},\\
    &\tau \in \{nT^0/2\}_{n=2,4,6,\ldots}.
\end{aligned}
$$
If $\tau \notin \{nT^0/2\}_{n\in \mathbb{N}_{\geq 0}}$, according to \eqref{eq:taus0}, \eqref{eq:C0}, and the continuity of $\boldsymbol{S}_{{\rm slow},\tau}^0 (-2,-2/3)$ in terms of $\tau$, one can calculate the limit
$$
\begin{aligned}
  \lim_{ \substack{(\bar{x},\bar{y})\rightarrow (1,-2/3)\\ (\bar{x},\bar{y})\in W_+} } \bar{f} ( \boldsymbol{S}_{{\rm slow},\tau-\tau_{\rm s}(\bar{x})}^0 (-2,-2/3)) \quad
  \\= \bar{f} ( \boldsymbol{S}_{{\rm slow},\tau}^0 (-2,-2/3)),
\end{aligned}
$$
which leads to \eqref{eq:Usingularm}. 
If $\tau \in \{nT^0/2\}_{n=1,3,5,\ldots}$, from \eqref{eq:sol3}
\begin{equation}
\label{eq:limit}
  \lim_{ \substack{(\bar{x},\bar{y})\rightarrow (1,-2/3)\\ (\bar{x},\bar{y})\in W_+} }  \boldsymbol{S}_{{\rm slow},\tau-\tau_{\rm s}(\bar{x})}^0 (-2,-2/3) = (-1,2/3),
\end{equation}
holds. 
Thus, one can calculate
$$
\begin{aligned}
  \lim_{ \substack{(\bar{x},\bar{y})\rightarrow (1,-2/3)\\ (\bar{x},\bar{y})\in W_+} } \bar{f} ( \boldsymbol{S}_{{\rm slow},\tau-\tau_{\rm s}(\bar{x})}^0 (-2,-2/3)) 
  &= \lim_{ \substack{(\bar{x},\bar{y})\rightarrow (-1,2/3)\\ (\bar{x},\bar{y})\in W_-} } \bar{f}(\bar{x},\bar{y})  \quad (\,\because~\text{\eqref{eq:limit}}\,) \\ 
  &= \bar{f} (2,2/3)  \quad (\,\because~\text{\eqref{eq:singularp}}\,) \\
  &= \bar{f} ( \boldsymbol{S}_{{\rm slow},\tau}^0 (-2,-2/3)),
\end{aligned}
$$
which leads to \eqref{eq:Usingularm}.
It can be shown similarly if $\tau \in \{nT^0/2\}_{n=2,4,6,\ldots}$. 
Therefore, \eqref{eq:Usingularm} holds for all $\tau> 0$ and $\bar{f}\in \bar{\mathcal{F}}$. (C2) is similarly proved. 
Therefore, we see (C4)-(C6) $\Longrightarrow$ (C1)-(C3) and \eqref{eq:slow_invariance}.
\end{proof}

\section{Proof of Lemma 2.}
\begin{proof}
First, we will prove $ \{{\rm i} n \omega^0 \}_{n\in \mathbb{Z}} = \sigma (U_\tau^0)$. 
To show $\{{\rm i} n \omega^0 \}_{n\in \mathbb{Z}} \subset \sigma (U_\tau^0)$, we use the notation $\left\{ \bar{\phi}_{{\rm i}\omega}^0 \right\}^n (\bar{x},\bar{y}):= \left\{ \bar{\phi}_{{\rm i}\omega}^0 (\bar{x},\bar{y}) \right\}^n $. 
For any $n\in \mathbb{Z}$, $\left\{ \bar{\phi}_{{\rm i}\omega}^0 \right\}^n \in \bar{\mathcal{F}}$, and we see that
\begin{equation}
\label{eq:multi_KEF}
\begin{aligned}
    U_{{\rm slow},\tau}^0 \left\{ \bar{\phi}_{{\rm i}\omega}^0 \right\}^n (\bar{x},\bar{y})&= \left\{ \bar{\phi}_{{\rm i}\omega}^0 (\boldsymbol{S}_{{\rm slow},\tau}^0 (\bar{x},\bar{y})) \right\}^n \\ 
    &= \left\{ U_{{\rm slow},\tau}^0 \bar{\phi}_{{\rm i}\omega}^0 (\bar{x},\bar{y}) \right\}^n \\
    &= \left\{ {\rm e}^{{\rm i} \omega^0 \tau} \bar{\phi}_{{\rm i}\omega}^0 (\bar{x},\bar{y}) \right\}^n \\&  = {\rm e}^{{\rm i}n \omega^0 \tau} \left\{ \bar{\phi}_{{\rm i}\omega}^0  \right\}^n (\bar{x},\bar{y}),
\end{aligned}
\end{equation}
holds for $\forall (\bar{x},\bar{y})\in W_\mp,~ \forall \tau \geq 0$; so $\{{\rm i} n \omega^0 \}_{n\in \mathbb{Z}} \subset \sigma (U_\tau^0)$.

To show $\{{\rm i} n \omega^0 \}_{n\in \mathbb{Z}} \supset \sigma (U_\tau^0)$,
we will show $\forall \lambda \in \mathbb{C}\backslash \{{\rm i} n \omega^0 \}_{n\in \mathbb{Z}}$ does not satisfy \eref{eq:slowKEF} which results in $\mathbb{C}\backslash \{{\rm i} n \omega^0 \}_{n\in \mathbb{Z}} \subset \mathbb{C}\backslash \sigma (U_\tau^0) $. 
Assume that there exists a $\lambda \in \mathbb{C}\backslash \{{\rm i} n \omega^0 \}_{n\in \mathbb{Z}}$ and $\bar{\phi}_\lambda \in \bar{\mathcal{F}}\backslash\{0\}$ satisfying \eref{eq:slowKEF}. 
The constrained flow $\boldsymbol{S}_{{\rm slow},\tau}^0$ defined in Appendix A satisfies that
$\forall (\bar{x},\bar{y}) \in W_\mp$, $\exists \tau>0$;
\begin{equation}
\nonumber    
\boldsymbol{S}_{\tau + T^0}^0 (\bar{x},\bar{y}) = \boldsymbol{S}_\tau^0 (\bar{x},\bar{y})
\end{equation}
where $T^0 = 2\pi/\omega^0$ is the period (which implies that all of the solutions become entirely periodic in a finite time or, in other words, that the solution is not unique). 
Thus, $\forall \bar{f}\in\bar{\mathcal{F}}$ satisfies
$$
\begin{aligned}
U_{{\rm slow},\tau + T^0}^0 \bar{f} (\bar{x},\bar{y}) &= \bar{f}\circ \boldsymbol{S}_{{\rm slow},\tau + T^0}^0 (\bar{x},\bar{y}) = \bar{f}\circ \boldsymbol{S}_{{\rm slow},\tau}^0 (\bar{x},\bar{y})\\ &= U_{{\rm slow},\tau}^0 \bar{f} (\bar{x},\bar{y}),
\end{aligned}
$$
Substituting $\bar{f}=\bar{\phi}_\lambda$ into the above, from the assumption of \eref{eq:slowKEF}, we obtain
$$
\begin{aligned}
    U_{{\rm slow},\tau}^0 \bar{\phi} (\bar{x},\bar{y}) &= {\rm e}^{\lambda T_0} U_{{\rm slow},\tau}^0 \bar{\phi} (\bar{x},\bar{y}) \\ &= {\rm e}^{\lambda \tau} \bar{\phi} (\bar{x},\bar{y}).
\end{aligned}
$$
Due to ${\rm e}^{\lambda T_0} \neq 1$, $U_{{\rm slow},\tau}^0 \bar{\phi} (\bar{x},\bar{y}) = \bar{\phi} (\bar{x},\bar{y}) = 0$ for $\forall (\bar{x},\bar{y})\in W_\mp$ holds and then we obtain $\bar{\phi}_\lambda = 0$ which contradicts the assumption $\bar{\phi}_\lambda \in \bar{\mathcal{F}}\backslash\{0\}$. 
Thus, we obtain $ \{{\rm i} n \omega^0 \}_{n\in \mathbb{Z}} \supset \sigma (U_\tau^0)$ and see $ \{{\rm i} n \omega^0 \}_{n\in \mathbb{Z}} = \sigma (U_\tau^0)$. 

Second, we prove the smoothness of the eigenfunctions. 
From \eqref{eq:slowKEFdef}, it is obvious that $\bar{\phi}_{{\rm i}\omega}^0|_{W_+}$ and $\bar{\phi}_{{\rm i}\omega}^0|_{W_-}$ are smooth. 
From \eqref{eq:multi_KEF}, it is also obvious that $\bar{\phi}_{{\rm i}n\omega}^0|_{W_+}$ and $\bar{\phi}_{{\rm i}n\omega}^0|_{W_-}$ are smooth. 
That completes the proof.

\end{proof} 

\bibliographystyle{unsrt}  
\bibliography{references.bib}

\begin{thebibliography}{10}

\bibitem{analysis_of_fluid_flowd_KO}
I.~Mezi{\'c}.
\newblock Analysis of fluid flows via spectral properties of the {K}oopman operator.
\newblock {\em Annual Review of Fluid Mechanics}, 45:357--378, 2013.

\bibitem{brunton2021modern}
Steven~L Brunton, Marko Budi{\v{s}}i{\'c}, Eurika Kaiser, and J~Nathan Kutz.
\newblock Modern {K}oopman theory for dynamical systems.
\newblock {\em SIAM Review}, 64(2):229--340, 2022.

\bibitem{otto2021koopman}
Samuel~E Otto and Clarence~W Rowley.
\newblock Koopman operators for estimation and control of dynamical systems.
\newblock {\em Annual Review of Control, Robotics, and Autonomous Systems}, 4:59--87, 2021.

\bibitem{lasota1998chaos}
Andrzej Lasota and Michael~C Mackey.
\newblock {\em {Chaos, Fractals, and Noise: Stochastic Aspects of Dynamics}}.
\newblock Springer-Verlag, New York, 1994.

\bibitem{arnold1968ergodic}
Vladimir~Igorevich Arnold and Andr{\'e} Avez.
\newblock Ergodic problems of classical mechanics.
\newblock {\em W.A.Benjamin, New York, Amsterdam}, 1968.

\bibitem{mezic2005spectral}
Igor Mezi{\'c}.
\newblock Spectral properties of dynamical systems, model reduction and decompositions.
\newblock {\em Nonlinear Dynamics}, 41:309--325, 2005.

\bibitem{rowley2009spectral}
Clarence~W Rowley, Igor Mezi{\'c}, Shervin Bagheri, Philipp Schlatter, and Dan~S Henningson.
\newblock Spectral analysis of nonlinear flows.
\newblock {\em {Journal of Fluid Mechanics}}, 641:115--127, 2009.

\bibitem{koopman-springer}
Alexandre Mauroy, Igor Mezi{\'c}, and Yoshihiko Susuki, editors.
\newblock {\em {The Koopman Operator in Systems and Control: Concepts, Methodologies, and Applications}}.
\newblock Springer, 2020.

\bibitem{use_of_Fourier_averages}
A.~Mauroy and I.~Mezi{\'c}.
\newblock On the use of {F}ourier averages to compute the global isochrons of (quasi) periodic dynamics.
\newblock {\em {CHAOS}: {An} {I}nterdisciplinary {J}ournal of {N}onlinear {S}cience}, 22(3):033112, 2012.

\bibitem{isostable_reduction_of_periodic_orbits}
D.~Wilson and J.~Moehlis.
\newblock Isostable reduction of periodic orbits.
\newblock {\em Physical Review E}, 94(5):052213, 2016.

\bibitem{phase-amplitude_reduction_nakao}
Sho Shirasaka, Wataru Kurebayashi, and Hiroya Nakao.
\newblock {P}hase-amplitude reduction of transient dynamics far from attractors for limit-cycling systems.
\newblock {\em {CHAOS: An Interdisciplinary Journal of Nonlinear Science}}, 27(2):023119, 2017.

\bibitem{global_computation_of_phase-amplitude_reduction}
A.~Mauroy and I.~Mezi{\'c}.
\newblock Global computation of phase-amplitude reduction for limit-cycle dynamics.
\newblock {\em {CHAOS: An Interdisciplinary Journal of Nonlinear Science}}, 28(7):073108, 2018.

\bibitem{mishchenko2013differential}
E~Mishchenko.
\newblock {\em {Differential Equations with Small Parameters and Relaxation Oscillations}}.
\newblock Springer Science \& Business Media, 1980.

\bibitem{izhikevich2000phase}
Eugene~M Izhikevich.
\newblock Phase equations for relaxation oscillators.
\newblock {\em SIAM Journal on Applied Mathematics}, 60(5):1789--1804, 2000.

\bibitem{desroches2012mixed}
Mathieu Desroches, John Guckenheimer, Bernd Krauskopf, Christian Kuehn, Hinke~M Osinga, and Martin Wechselberger.
\newblock Mixed-mode oscillations with multiple time scales.
\newblock {\em SIAM Review}, 54(2):211--288, 2012.

\bibitem{kuehn2015multiple}
Christian Kuehn.
\newblock {\em {Multiple Time Scale Dynamics}}, volume 191.
\newblock Springer, 2015.

\bibitem{isostable_isochron_and_koopmanSpectrum_fixedpoint}
A.~Mauroy, I.~Mezi{\'c}, and J.~Moehlis.
\newblock Isostables, isochrons, and {K}oopman spectrum for the action--angle representation of stable fixed point dynamics.
\newblock {\em Physica D: Nonlinear Phenomena}, 261:19--30, 2013.

\bibitem{eldering2018global}
Jaap Eldering, Matthew Kvalheim, and Shai Revzen.
\newblock Global linearization and fiber bundle structure of invariant manifolds.
\newblock {\em Nonlinearity}, 31(9):4202, 2018.

\bibitem{fenichel1979geometric}
Neil Fenichel.
\newblock Geometric singular perturbation theory for ordinary differential equations.
\newblock {\em {Journal of Differential Equations}}, 31(1):53--98, 1979.

\bibitem{singular_perturbation_geometric}
Christopher~KRT Jones.
\newblock Geometric singular perturbation theory.
\newblock {\em {Dynamical Systems: Lectures Given at the 2nd Session of the Centro Internazionale Matematico Estivo (CIME) held in Montecatini Terme, Italy, June 13--22, 1994}}, pages 44--118, 1995.

\bibitem{linearization_in_the_large_of_nonlinear_systems_and_KoopmanOperatorSpectrum}
Y.~Lan and I.~Mezi{\'c}.
\newblock Linearization in the large of nonlinear systems and {K}oopman operator spectrum.
\newblock {\em Physica D: Nonlinear Phenomena}, 242(1):42--53, 2013.

\bibitem{osinga2010continuation}
Hinke~M Osinga and Jeff Moehlis.
\newblock Continuation-based computation of global isochrons.
\newblock {\em SIAM Journal on Applied Dynamical Systems}, 9(4):1201--1228, 2010.

\bibitem{sherwood2010dissecting}
William~Erik Sherwood and John Guckenheimer.
\newblock Dissecting the phase response of a model bursting neuron.
\newblock {\em SIAM Journal on Applied Dynamical Systems}, 9(3):659--703, 2010.

\bibitem{mauroy2014global}
Alexandre Mauroy, Blane Rhoads, Jeff Moehlis, and Igor Mezic.
\newblock Global isochrons and phase sensitivity of bursting neurons.
\newblock {\em SIAM Journal on Applied Dynamical Systems}, 13(1):306--338, 2014.

\bibitem{susuki2021koopman}
Yoshihiko Susuki.
\newblock On {K}oopman operator framework for semi-explicit differential-algebraic equations.
\newblock {\em IFAC-PapersOnLine}, 54(14):341--345, 2021.

\bibitem{takens1976constrained}
Floris Takens.
\newblock Constrained equations; a study of implicit differential equations and their discontinuous solutions.
\newblock In {\em {Structural Stability, the Theory of Catastrophes, and Applications in the Sciences: Proceedings of the Conference Held at Battelle Seattle Research Center 1975}}, pages 143--234. Springer, 1976.

\bibitem{sastry1981jump}
Shankar Sastry and C~Desoer.
\newblock Jump behavior of circuits and systems.
\newblock {\em IEEE Transactions on Circuits and Systems}, 28(12):1109--1124, 1981.

\bibitem{spectrum-of-the-koopman-operator}
I.~Mezi{\'c}.
\newblock Spectrum of the {K}oopman operator, spectral expansions in functional spaces, and state-space geometry.
\newblock {\em Journal of Nonlinear Science}, 30(5):2091--2145, 2020.

\bibitem{Existence_Kvalheim}
M.~D. Kvalheim and S.~Revzen.
\newblock Existence and uniqueness of global {K}oopman eigenfunctions for stable fixed points and periodic orbits.
\newblock {\em Physica D: Nonlinear Phenomena}, 425:132959, 2021.

\bibitem{krupa2001relaxation}
Martin Krupa and Peter Szmolyan.
\newblock Relaxation oscillation and canard explosion.
\newblock {\em Journal of Differential Equations}, 174(2):312--368, 2001.

\bibitem{gaspard1995spectral}
Pierre Gaspard, Gr{\'e}goire Nicolis, Astero Provata, and S13835391995PhRvE Tasaki.
\newblock Spectral signature of the pitchfork bifurcation: Liouville equation approach.
\newblock {\em Physical Review E}, 51(1):74, 1995.

\bibitem{shirasaka2017phase}
Sho Shirasaka, Wataru Kurebayashi, and Hiroya Nakao.
\newblock Phase reduction theory for hybrid nonlinear oscillators.
\newblock {\em Physical Review E}, 95(1):012212, 2017.

\bibitem{govindarajan2016operator}
Nithin Govindarajan, Hassan Arbabi, Louis Van~Blargian, Timothy Matchen, Emma Tegling, and I.~Mezi\'{c}.
\newblock An operator-theoretic viewpoint to non-smooth dynamical systems: Koopman analysis of a hybrid pendulum.
\newblock In {\em 2016 IEEE 55th Conference on Decision and Control (CDC)}, pages 6477--6484. IEEE, 2016.

\end{thebibliography}

\end{document}